\documentclass[10pt]{amsart}
\usepackage{amsmath}
\usepackage{amssymb}
\usepackage{amsfonts}
\usepackage{amsthm}
\usepackage{amsxtra}
\usepackage{portland}
\usepackage{rotating}
\usepackage{nicefrac}
\usepackage{float}
\usepackage[all,web,arc,poly,dvips]{xy}
\usepackage{epic,eepic}
\usepackage{epsfig}
\usepackage[dvips]{color}
\usepackage{array}
\usepackage{pifont}
\usepackage{multirow}
\usepackage{latexsym}
\usepackage{graphics}
\usepackage{fullpage}
\usepackage{hyperref}

\theoremstyle{plain}

\newtheorem{corollary}{Corollary}[section]

\newtheorem{proposition}{Proposition}[section]
\theoremstyle{definition}
\newtheorem{definition}{Definition}[section]

\theoremstyle{remark}
\newtheorem{remark}{Remark}[section]

\DeclareMathAlphabet{\mathpzc}{OT1}{pzc}{m}{it}
\newcommand{\K}{\mathpzc k}
\newcommand{\C}{\mathbb C}

\newcommand{\Z}{\mathbb Z}

\newcommand{\T}{\mathbb T}

\newcommand{\half}{
        {\lower0.00ex\hbox{\raise.6ex\hbox{\the\scriptfont0 1}
                           \kern-.5em\slash\kern-.1em\lower.45ex
                                     \hbox{\the\scriptfont0 2}}}}
\newcommand{\quarter}{
        {\lower0.00ex\hbox{\raise.6ex\hbox{\the\scriptfont0 1}
                           \kern-.5em\slash\kern-.1em\lower.45ex
                                     \hbox{\the\scriptfont0 4}}}}
\newcommand{\tquarter}{
        {\lower0.00ex\hbox{\raise.6ex\hbox{\the\scriptfont0 3}
                           \kern-.5em\slash\kern-.1em\lower.45ex
                                     \hbox{\the\scriptfont0 4}}}}
\newcommand{\eighth}{
        {\lower0.00ex\hbox{\raise.6ex\hbox{\the\scriptfont0 1}
                           \kern-.5em\slash\kern-.1em\lower.45ex
                                     \${\rm P}_{\rm III^{\prime}}\;$hbox{\the\scriptfont0 8}}}}
\newcommand{\othird}{
        {\lower0.00ex\hbox{\raise.6ex\hbox{\the\scriptfont0 1}
                           \kern-.5em\slash\kern-.1em\lower.45ex
                                     \hbox{\the\scriptfont0 3}}}}

\def\cn{{\rm cn}}
\def\dn{{\rm dn}}
\def\sn{{\rm sn}}

\def\eK{{\rm K}}
\def\eE{{\rm E}}

\def\eZeta{\mathcal{Z}}

\begin{document}

\title[...]{Fredholm {D}eterminant evaluations of the Ising Model diagonal correlations and their $ \lambda $ generalisation}

\author{N.S.~Witte \and P.J.~Forrester}
\address{Department of Mathematics and Statistics,
University of Melbourne,Victoria 3010, Australia}
\email{\tt n.witte@ms.unimelb.edu.au}

\begin{abstract}
The diagonal spin-spin correlations of the square lattice Ising model, originally expressed as Toeplitz 
determinants, are given by two distinct Fredholm
determinants - one with an integral operator having an Appell function kernel and another with a summation
operator having a Gauss hypergeometric function kernel. Either determinant allows for a Neumann expansion possessing
a natural $ \lambda $-parameter generalisation and we prove that both expansions are in fact equal, implying
a continuous and a discrete representation of the form factors. Our proof employs an extension of the classic
study by Geronimo and Case \cite{GC_1979}, applying scattering theory to orthogonal polynomial systems on the
unit circle, to the bi-orthogonal situation.
\end{abstract}

\maketitle

\section{Form Factor expansions of the Diagonal Correlations}
\label{FFE}
\setcounter{equation}{0}
The two-dimensional Ising model and the classical random matrix ensembles are examples of model
systems with wide ranging applicability. Taken literally, the two-dimensional Ising model
specifies a classical magnetic system from the microscopic interactions of two-state
spin variables. These variables are confined to a particular two-dimensional lattice,
and the interactions are confined to neighbouring lattice sites. In particular, in the
case of the square lattice, each site $ (i,j) $ carries a spin variable 
$ \sigma_{i,j} \in \{-1,1\} $ with coupling between nearest neighbours in the 
horizontal and vertical directions. The joint probability density function for a 
particular configuration $ \{\sigma_{i,j}\}_{-N \leq i,j \leq N} $ of the states on a
$ (2N+1)\times(2N+1) $ lattice centred at the origin is given by
\begin{equation} 
   P_{2N+1}(\{\sigma_{i,j}\}) = \frac{1}{Z_{2N+1}}
  \exp\left( K_1\sum^{N-1}_{i=-N}\sum^{N}_{j=-N}\sigma_{i,j}\sigma_{i+1,j}+K_2\sum^{N}_{i=-N}\sum^{N-1}_{j=-N}\sigma_{i,j}\sigma_{i,j+1} \right) ,
\label{ising_config}
\end{equation}
where $ Z_{2N+1} $ is the normalisation such that summing (\ref{ising_config}) over all allowed $ \{\sigma_{i,j}\} $
gives unity. The parameters $ K_1 $ and $ K_2 $ are dimensionless couplings in the horizontal 
and vertical directions respectively. As a model of a magnet it would be natural to
generalise (\ref{ising_config}) to include a coupling to an external magnetic field 
which adds
\begin{equation*}
   h\sum^{N}_{i,j=-N}\sigma_{i,j}
\end{equation*}
to the exponent in (\ref{ising_config}) to give the joint probability density
$ P^{h}_{2N+1}(\{\sigma_{i,j}\}) $. In the absence of such a term one is said to
be considering the zero field Ising model.

It is a celebrated result due to Peierls that the zero field two-dimensional Ising model
exhibits a phase transition from a high temperature disordered phase characterised by
zero spontaneous magnetisation to a low temperature ordered phase with non-zero
magnetisation. Setting $ \K=\sinh 2K_1 \sinh 2K_2 $ the Kramers-Wannier duality argument
gives that the phase transition occurs at $ \K=1 $. The zero field magnetisation is
specified by
\begin{equation}
   \langle \sigma_{0,0} \rangle = 
  \lim_{h\to 0}\lim_{N\to \infty} \sum_{\{\sigma_{i,j}\}}\sigma_{0,0}P^{h}_{2N+1}(\{\sigma_{i,j}\}) .
\label{magnet}
\end{equation}
Moreover there is a closed form expression for $ \langle \sigma_{0,0} \rangle $ in the
low-temperature phase
\begin{equation}
   \langle \sigma_{0,0} \rangle = (1-\K^{-2})^{1/8} ,
\label{low_magnet}
\end{equation}
which was announced by Onsager in 1948 \cite{Onsager_1948} and proven by Yang in 1952 \cite{Yang_1952}.

Our interest is in the zero-field spin-spin correlation function between spins 
$ \sigma_{0,0} $ at the centre of the lattice, and the spin $ \sigma_{i,j} $ at
site $ (i,j) $. In the infinite lattice limit this is defined as
\begin{equation}
   \langle \sigma_{0,0}\sigma_{i,j} \rangle = 
  \lim_{N\to \infty} \sum_{\{\sigma_{i,j}\}}\sigma_{0,0}\sigma_{i,j}P_{2N+1}(\{\sigma_{i,j}\}) .
\label{spin_spin}
\end{equation}
Onsager and Kaufman knew of a Toeplitz determinant form for the case $ i=j=n $ (spin-spin correlation along the diagonal)
but never published their result. However a draft paper has recently come to light which obtains this result
and, according to Baxter \cite{Baxter_2011}, was almost certainly written by the pair. This reads \cite{McCW_1973}
\begin{equation}
   \langle \sigma_{0,0}\sigma_{n,n} \rangle = 
  \det [a_{i-j} ]_{1\leq i,j\leq n}, \quad 
   a_{p}=\int^{\pi}_{-\pi}\frac{d\theta}{2\pi}e^{-ip\theta}a(\theta) ,
\label{diagonal}
\end{equation}
where
\begin{equation}
   a(\theta) = \left[ \frac{1-\K^{-1}e^{-i\theta}}{1-\K^{-1}e^{i\theta}} \right]^{1/2} ,
\label{diagonal:a}
\end{equation}
and $ n\in \mathbb{N} $. This result (\ref{diagonal}) in fact includes (\ref{low_magnet}) as
a special case. To see this one notes that for $ n $ large
\begin{equation}
  \langle \sigma_{0,0}\sigma_{n,n} \rangle \sim \langle \sigma_{0,0} \rangle\langle \sigma_{n,n} \rangle =  \langle \sigma_{0,0} \rangle^{2} ,
\label{factor}
\end{equation}
upon the assumption that to leading order spins far apart are independent of each other.
On the other hand for Toeplitz determinants with symbol $ a(\theta) $ such that $ \log a(\theta) $
is Laurent expandable in $ e^{i\theta} $ 
\begin{equation}
   \log a(\theta) = \sum^{\infty}_{p=-\infty} c_{p}e^{ip\theta} ,
\label{LogWgt}
\end{equation}
and with the further requirement that $ \sum^{\infty}_{p=-\infty} |p||c_{p}|^2 < \infty $, 
the Szeg\H{o} theorem \cite{Jo_1988} gives
\begin{equation}
  \det [a_{i-j} ]_{i,j=1,\ldots,n} \mathop{\sim}\limits_{n \to \infty}
   \exp\left( nc_{0}+\sum^{\infty}_{p=1}pc_{p}c_{-p}+{\rm o}(1) \right) .
\label{szego}
\end{equation}
For $ a(\theta) $ given by (\ref{diagonal:a}) and $ \K>1 $ (low temperature phase)
\begin{equation}
   c_0 = 0, \qquad c_{p} = \frac{1}{2p}\K^{-|p|}, \quad p\in \Z,\, p\neq 0 .
\label{LogFCoeff}
\end{equation}
The Szeg\H{o} theorem then gives
\begin{equation}
   \langle \sigma_{0,0}\sigma_{n,n} \rangle \sim (1-\K^{-2})^{1/4} 
\label{onsager44}
\end{equation}
which when substituted in (\ref{factor}) implies (\ref{low_magnet}).

Intimately related to the asymptotic expansion (\ref{szego}) is the so-called form factor
expansion
\begin{equation}
   \langle \sigma_{0,0}\sigma_{n,n} \rangle = (1-\K^{-2})^{1/4}\left( 1+\sum^{\infty}_{p=1} f^{(2p)}_{n,n} \right) 
\label{ffexp}
\end{equation}
and its $ \lambda $ generalisation
\begin{equation}
   C_{-}(n,n;\lambda) = (1-\K^{-2})^{1/4}\left( 1+\sum^{\infty}_{p=1} \lambda^{2p} f^{(2p)}_{n,n} \right)  .
\label{ffexp2}
\end{equation}
Expansions of the type (\ref{ffexp}) were initiated by Wu \cite{Wu_1966} in the study of the 
large $ n $ form of (\ref{spin_spin}) with $ i=0, j=n $ so the spins are along the same row.
Each term in (\ref{ffexp}) is to be a higher order correction in
the large $ n $ limit; and to leading order in $ \K^{-2} $ is proportional to $ \K^{-2p(n+p)} $.
The coefficients $ f^{(2p)}_{n,n} $ furthermore have the interpretation of $2p$ quasi-particle
contributions to the two-point correlation function \cite{WMTB_1976}

An explicit $2p$-dimensional multiple integral formula for $ f^{(2p)}_{n,n} $ follows as 
the special case of the explicit form of the form factor expansion of
$ \langle \sigma_{0,0}\sigma_{n,n} \rangle $ given in \cite{WMTB_1976}. However this
is complicated, and in fact starting directly with (\ref{diagonal}) the simpler expression
\begin{multline}
   f^{(2p)}_{n,n} = \frac{t^{p(n+p)}}{(p!)^2\pi^{2p}}\int^1_0dx_1 \cdots \int^1_0dx_{2p}
  \\
  \times \prod^{2p}_{k=1}x^{n}_{k}\prod^{p}_{j=1}\left[ \frac{(1-tx_{2j})(x^{-1}_{2j}-1)}{(1-tx_{2j-1})(x^{-1}_{2j-1}-1)} \right]^{1/2}
  \\
  \times \prod^{p}_{j=1}\prod^{p}_{k=1}(1-tx_{2k-1}x_{2j})^{-2}\prod_{1\leq j<k \leq p}(x_{2j-1}-x_{2k-1})^2(x_{2j}-x_{2k})^2 ,
\label{ff_Ieven}
\end{multline}
where $ t=\K^{-2}<1 $ and $ n\geq 0 $, has been derived originally by \cite{Bu_2001}, \cite{BL_2004} and
independently by \cite{LM_2007}, \cite{BHMMZ_2007}.
In the high temperature regime the correlations are given by
\begin{equation}
   \langle \sigma_{0,0}\sigma_{n,n} \rangle = (1-t)^{1/4}\sum^{\infty}_{p=0} f^{(2p+1)}_{n,n}
\label{ffhigh}
\end{equation}
and
\begin{multline}
   f^{(2p+1)}_{n,n} = \frac{t^{n(p+1/2)+p(p+1)}}{p!(p+1)!\pi^{2p+1}}\int^1_0dx_1 \cdots \int^1_0dx_{2p+1}
  \\
  \times \prod^{2p+1}_{k=1}x^{n}_{k}\prod^{p+1}_{j=1}\frac{1}{x_{2j-1}}(1-tx_{2j-1})^{-1/2}(x^{-1}_{2j-1}-1)^{-1/2}
                                    \prod^p_{j=1}x_{2j}(1-tx_{2j})^{1/2}(x^{-1}_{2j}-1)^{1/2}
  \\
  \times \prod^{p+1}_{j=1}\prod^{p}_{k=1}(1-tx_{2j-1}x_{2k})^{-2}\prod_{1\leq j<k \leq p+1}(x_{2j-1}-x_{2k-1})^2\prod_{1\leq j<k\leq p}(x_{2j}-x_{2k})^2 ,
\label{ff_Iodd}
\end{multline}
and in this case $ t=\K^2<1 $ and $ n\geq 0 $.

It is our aim to develop the analytic properties of (\ref{ffexp}-\ref{ff_Iodd}), building on
known results from \cite{FW_2004b}, \cite{BHMMOZ_2007}, \cite{MG_2010}. In particular a theory
will be developed which leads to two distinct Fredholm determinant formulas for the 
$ \lambda $ generalised form factor expansion (\ref{ffexp2}).
We begin in Section \ref{CFD} by evaluating the first two terms in the power series
expansion of (\ref{ff_Ieven}) and writing (\ref{ffexp2}) as a Fredholm determinant of a particular integral 
operator supported on $ (0,1) $. In Section \ref{DFD}
we express (\ref{ffexp}) as a discrete Fredholm determinant, making use of the 
Borodin-Okounkov identity \cite{GC_1979}, \cite{BO_2000a} from the theory of determinantal
point processes \cite{rmt_Fo} and we conjecture that with the kernel multiplied by $ \lambda^2 $
(\ref{ffexp2}) results. The rest of our study is devoted to proving this conjecture and thus the equality of
the two Fredholm determinants, although our working is more general. 
Thus our first step in this undertaking is to introduce the
bi-orthogonal polynomial system attached to a general weight $ a(\zeta) $ in Section
\ref{BOPS}. Next we indicate in Section \ref{ST} how the classic results of Geronimo and Case \cite{GC_1979} can
be generalised to non-hermitian Toeplitz determinants.
Finally in Section \ref{IM} we pose a $ \lambda $ extension of a key functional equation arising in 
Section \ref{BOPS} and deduce some consequences of this, in particular we describe how this is connected
with the pioneering work of Wu \cite{Wu_1966}. It is through this connection that we
can establish the desired equality.

\section{Fredholm Determinant with Appell Kernel} 
\label{CFD}
\setcounter{equation}{0}
\subsection{$ t\to 0 $ Expansions of form factors}
We see from (\ref{ff_Ieven}, \ref{ff_Iodd}) that $ f^{(2p)}_{n,n}, f^{(2p+1)}_{n,n} $ have the small $ t $ expansions
\begin{equation}
  f^{(2p)}_{n,n} = t^{p(n+p)}\sum^{\infty}_{k=0}c^{<}_{k}t^k , \qquad
  f^{(2p+1)}_{n,n} = t^{n(p+1/2)+p(p+1)}\sum^{\infty}_{k=0}c^{>}_{k}t^k ,
\label{ffsmallt}
\end{equation}
where $ c^{<,>}_k=c^{<,>}_k(n,p) $ and $ t=\K^{\mp 2} $ respectively. 
Using knowledge of Selberg integral theory (see e.g. \cite{rmt_Fo}) the first two coefficients of 
(\ref{ffsmallt}) can be readily computed.
\begin{proposition}
For $ t=\K^{-2}<1 $ we have
\begin{equation}
   c^{<}_0 = \frac{1}{(p!)^2\pi^{2p}}\frac{\Gamma(n+p+\tfrac{1}{2})\Gamma(p+\tfrac{1}{2})}{\Gamma(n+\tfrac{1}{2})\Gamma(\tfrac{1}{2})}
   \left( \prod^{p-1}_{j=0}\frac{\Gamma(n+j+\tfrac{1}{2})\Gamma(j+\tfrac{1}{2})\Gamma(j+2)}{\Gamma(n+p+j+1)} \right)^2 ,
\label{ffzero:a}
\end{equation}
and
\begin{equation}
   c^{<}_1 = c^{<}_0 \frac{p(n+p)}{2(n+2p)^2}[4p(n+p)+1] ,
\label{ffone:a}
\end{equation}
while for $ t=\K^2<1 $ we have
\begin{equation}
   c^{>}_0 = \frac{1}{p!\pi^{2p+1}}\frac{\Gamma(n+\tfrac{1}{2})\Gamma(\tfrac{1}{2})}{\Gamma(n+p+1)}
   \left( \prod^{p-1}_{j=0}\frac{\Gamma(n+j+\tfrac{3}{2})\Gamma(j+\tfrac{3}{2})\Gamma(j+2)}{\Gamma(n+p+j+2)} \right)^2 ,
\label{ffzero:b}
\end{equation}
and
\begin{equation}
   c^{>}_1 = c^{>}_0 \frac{(n+p+\tfrac{1}{2})}{2(n+2p+1)^2}[4p(p+1)(n+p+\tfrac{1}{2})+n+2p+1] .
\label{ffone:b}
\end{equation}
\end{proposition}
\begin{proof}
We read off from (\ref{ff_Ieven}) that
\begin{multline}
   c_{0} = \frac{1}{(p!)^2\pi^{2p}}\int^1_0dx_1 \cdots \int^1_0dx_{2p}
  \prod^{2p}_{k=1}x^{n}_{k}\prod^{p}_{j=1}\left[ \frac{(x^{-1}_{2j}-1)}{(x^{-1}_{2j-1}-1)} \right]^{1/2}
  \prod_{1\leq j<k \leq p}(x_{2j-1}-x_{2k-1})^2(x_{2j}-x_{2k})^2 ,
  \\
  = \frac{1}{(p!)^2\pi^{2p}}
    \left( \int^1_0dx_2\int^1_0dx_4 \cdots \int^1_0dx_{2p}
           \prod^{p}_{k=1}x^{n-1/2}_{2k}(1-x_{2k})^{1/2}
           \prod_{1\leq j<k \leq p}(x_{2j}-x_{2k})^2 \right)
  \\ \times
    \left( \int^1_0dx_1\int^1_0dx_{3} \cdots \int^1_0dx_{2p-1}
           \prod^{p}_{k=1}x^{n+1/2}_{2k-1}(1-x_{2k-1})^{-1/2}
           \prod_{1\leq j<k \leq p}(x_{2j-1}-x_{2k-1})^2 \right) .
\label{selberg1}
\end{multline}
We recognise both integrals in (\ref{selberg1}) as special cases of the Selberg integral
\begin{equation}
    S_{N}(\lambda_1,\lambda_2,\lambda) := 
    \int^1_0dt_1 \cdots \int^1_0dt_N \prod^N_{l=1}t_{l}^{\lambda_1}(1-t_{l})^{\lambda_2}\prod_{1\leq j<k\leq N}|t_{j}-t_{k}|^{2\lambda} ,
\label{selberg}
\end{equation}
which, according to Selberg \cite{Se_1944}, has the gamma function evaluation
\begin{equation*}
  S_{N}(\lambda_1,\lambda_2,\lambda) =
  \prod^{N-1}_{j=0}\frac{\Gamma(\lambda_1+1+j\lambda)\Gamma(\lambda_2+1+j\lambda)\Gamma(1+(j+1)\lambda)}
                        {\Gamma(\lambda_1+\lambda_2+2+(N+j-1)\lambda)\Gamma(1+\lambda)} . 
\end{equation*}
The result (\ref{ffzero:a}) follows.
In relation to (\ref{ffone:a}) we read off from (\ref{ff_Ieven}) that
\begin{equation*}
    c_{1} = \frac{1}{(p!)^2\pi^{2p}}\left( -\tfrac{1}{2}A_1+\tfrac{1}{2}A_2+2A_3 \right) ,
\end{equation*}
where, with $ S_{p}(\lambda_1,\lambda_2,\lambda)[f] $ denoting the Selberg integral (\ref{selberg})
with an additional factor $ f $ in the integrand
\begin{gather*}
   A_1 = S_{p}(n-\tfrac{1}{2},\tfrac{1}{2},1)\Big[\sum^p_{j=1}t_j\Big]S_{p}(n+\tfrac{1}{2},-\tfrac{1}{2},1) ,
  \\
   A_2 = S_{p}(n-\tfrac{1}{2},\tfrac{1}{2},1)S_{p}(n+\tfrac{1}{2},-\tfrac{1}{2},1)\Big[\sum^p_{j=1}t_j\Big] ,
  \\
   A_3 = S_{p}(n-\tfrac{1}{2},\tfrac{1}{2},1)\Big[\sum^p_{j=1}t_j\Big]S_{p}(n+\tfrac{1}{2},-\tfrac{1}{2},1)\Big[\sum^p_{j=1}t_j\Big] .
\end{gather*}
But according to a result of Aomoto \cite{Ao_1987}
\begin{equation*}
   \frac{1}{S_{p}(\lambda_1,\lambda_2,\lambda)}S_{p}(\lambda_1,\lambda_2,\lambda)\Big[\sum^p_{j=1}t_j\Big]
   = p\frac{\lambda_1+(p-1)\lambda+1}{\lambda_1+\lambda_2+2(p-1)\lambda+2} ,
\end{equation*}
which gives (\ref{ffone:a}). The corresponding results for the high temperature phase are found
by identical means.
\end{proof}

\subsection{Fredholm determinants with continuous kernel}
Our next result is to show that $ \sum^{\infty}_{p=0}f^{(2p)}_{n,n} $ with each term as given by 
(\ref{ff_Ieven}) can be expressed as a Fredholm determinant with kernel $ K(x,y) $ of the so-called
integrable type \cite{HI_2002}
\begin{equation*}
   K(x,y) = \frac{\phi(x)\psi(y)-\phi(y)\psi(x)}{x-y} .
\end{equation*}
\begin{proposition}
The low temperature diagonal spin-spin correlation of the square lattice Ising model
$ \langle \sigma_{0,0}\sigma_{n,n} \rangle $, $ |t|<1 $, $ n\geq 0 $ is, to within the prefactor (\ref{onsager44}),
the Fredholm determinant
\begin{equation}
   \det(\mathbb{I}-\mathbb{K}^{-}_{(0,1)}) = 1+\sum^{\infty}_{p=1}f^{(2p)}_{n,n} ,
\label{lowT_CFredholmDet}
\end{equation}
where $ \mathbb{K}^{-} $ is the integral operator with kernel $ K^{-}(x,y) $ of the explicit form
\begin{multline}
    K^{-}(x,y) = -\frac{\Gamma(n+\tfrac{1}{2})\Gamma(\tfrac{1}{2})}{2\pi^2(n+1)!}t^{n+1}(xy)^{n/2+1/4}
    \left[ (1-x)(1-y)(1-tx)(1-ty) \right]^{-1/4}
    \\ \times
    \frac{1}{x-y}\left[ xF_{1}(n+\tfrac{1}{2};-\tfrac{1}{2},1;n+2;t,tx)-yF_{1}(n+\tfrac{1}{2};-\tfrac{1}{2},1;n+2;t,ty) \right] ,
\label{lowT_Ckernel}
\end{multline}
and $ F_1(a;b,c;d;x,y) $ is the first Appell function \cite{DLMF}, \textsection 16.13.
Consequently we have
\begin{equation}
  f^{(2p)}_{n,n} = \frac{(-1)^p}{p!} \int^1_0du_1 \cdots \int^1_0du_p \det[ K^{-}(u_j,u_k) ]_{1\leq j,k\leq p} .
\label{NeumannLowT}
\end{equation}
\end{proposition}
\begin{proof}
Our proof parallels that taken by \cite{SGK_2008} in their treatment of a similar integral
for the impenetrable Bose gas.
For convenience we make the replacements $ x_{2j-1}\mapsto u_j, x_{2j}\mapsto v_j $ in (\ref{ff_Ieven})
and separate out the multiple $v$-integrations from those of the $u$-integrations thus
\begin{multline*}
   f^{(2p)}_{n,n} = \frac{t^{p(n+p)-2p^2}}{(p!)^2\pi^{2p}}\int^1_0du_1 \cdots \int^1_0du_p
   \prod^{p}_{j=1}u^{n-2p}_j\prod^{p}_{j=1}\left[(1-tu_j)(u^{-1}_j-1)\right]^{-1/2} \Delta^2_p(u)
   \\
   \times \int^1_0dv_1 \cdots \int^1_0dv_p \prod^{p}_{j=1}v^n_j\prod^{p}_{j=1}\left[(1-tv_j)(v^{-1}_j-1)\right]^{1/2}
          \Delta^2_p(v) \prod^{p}_{j,k=1}(v_j-t^{-1}u^{-1}_k)^{-2} ,
\end{multline*}
where $ \Delta_p(v) $ is the Vandermonde determinant or product of differences of the $ \{v_j\}^p_{j=1} $.
The inner $v$-integral can be written as a determinant using the Heine identity (see Eqs. (2.27) and (2.2.10)
of \cite{ops_Sz} and the original \cite{He_1961}) to give
\begin{equation}
   p! \det\left[ \int^1_0 dv\, v^{n+j+k-2}\frac{[(1-tv)(v^{-1}-1)]^{1/2}}{\prod^p_{l=1}(v-t^{-1}u^{-1}_l)^2} \right]_{1\leq j,k\leq p} .  
\label{hankel_aux}
\end{equation}
However this determinant form does not yield the simplest evaluation
so we employ the additional identity for general distinct indeterminates $ \{y_j\}^p_{j=1} $
\begin{equation}
   \det\left[ \int^1_0 dx\, x^{j+k-2}\frac{w(x)}{\prod^p_{l=1}(x-y_l)^2} \right]_{1\leq j,k\leq p}
  = \frac{1}{\Delta^2(y)} \det\left[ \int^1_0 dx\frac{w(x)}{(x-y_j)(x-y_k)} \right]_{1\leq j,k\leq p} .
\label{hankel_ID}
\end{equation}
Applying this to (\ref{hankel_aux}) we arrive at our simplest form, except for the evaluation
of the matrix elements. Using the integral representation of the first Appell function, i.e.
see \cite{DLMF}, Eq. (16.15.1) 
\begin{equation*}
   \frac{\Gamma(\alpha)\Gamma(\gamma-\alpha)}{\Gamma(\gamma)}F_{1}(\alpha;\beta,\beta';\gamma;x,y) 
  = \int^{1}_{0}du u^{\alpha-1}(1-u)^{\gamma-\alpha-1}(1-x u)^{-\beta}(1-y u)^{-\beta'}, \quad {\rm Re}(\alpha)>0, {\rm Re}(\gamma-\alpha)>0 ,
\end{equation*}
with
\begin{equation*}
   F_{1}(\alpha;\beta,\beta';\gamma;x,y) 
  = \sum^{\infty}_{m,n=0}\frac{(\alpha)_{m+n}(\beta)_{m}(\beta')_{n}}{(\gamma)_{m+n}}\frac{x^my^n}{m!n!}, \quad {\rm max}(|x|,|y|)<1 ,
\end{equation*}
we arrive at the result (\ref{lowT_CFredholmDet}) and (\ref{lowT_Ckernel}).
\end{proof}

There is an analogous result for the high temperature correlations although we have in
this case the added complication of an odd number of integrations and therefore an unpaired
variable.
\begin{proposition}
The high temperature diagonal spin-spin correlation of the square lattice Ising model
$ \langle \sigma_{0,0}\sigma_{n,n} \rangle $, $ |t|<1 $, $ n\geq 0 $ is, aside from the prefactor (\ref{onsager44}),
a Fredholm minor and its Neumann expansion coefficients are
\begin{equation}
   f^{(2p+1)}_{n,n} = \frac{1}{\pi p!}\int^1_0dv_1 \cdots \int^1_0dv_p
   \det\left[ \begin{matrix} K^{+}_0 & K^{+}_1(v_k) \cr K^{+}_1(v_j) & K^{+}_2(v_j,v_k) \end{matrix} \right]_{1\leq j,k\leq p} ,
   \quad p>0, \quad f^{(1)}_{n,n} = \frac{1}{\pi}K^{+}_0 ,
\label{highT_CFredholm}
\end{equation}
where the kernels are
\begin{multline}
    K^{+}_2(x,y) = \frac{\Gamma(n+\tfrac{1}{2})\Gamma(\tfrac{1}{2})}{\pi^2 n!}t^{n}(xy)^{n/2-3/4}
    \left[ (1-x)(1-y)(1-tx)(1-ty) \right]^{1/4}
    \\ \times
    \frac{1}{x-y}\left[ xF_{1}(n+\tfrac{1}{2};\tfrac{1}{2},1;n+1;t,tx)-yF_{1}(n+\tfrac{1}{2};\tfrac{1}{2},1;n+1;t,ty) \right] ,
\label{highT_Ckernel:a}
\end{multline}
\begin{equation}
    K^{+}_1(x) = -\frac{\Gamma(n+\tfrac{1}{2})\Gamma(\tfrac{1}{2})}{\pi n!}t^{n/2}x^{n/2-3/4} \left[ (1-x)(1-tx) \right]^{1/4}
    F_{1}(n+\tfrac{1}{2};\tfrac{1}{2},1;n+1;t,tx) ,
\label{highT_Ckernel:b}
\end{equation}
and 
\begin{equation}
    K^{+}_0 = \frac{\Gamma(n+\tfrac{1}{2})\Gamma(\tfrac{1}{2})}{n!}t^{n/2} {}_2F_{1}(n+\tfrac{1}{2},\tfrac{1}{2};n+1;t) .
\label{highT_Ckernel:c}
\end{equation}
\end{proposition}
\begin{proof}
For convenience we make the replacements $ x_{2j-1}\mapsto u_j, x_{2j}\mapsto v_j $ 
in (\ref{ff_Iodd}) and separate out the multiple $v$-integrations from those of the $u$-integrations 
making the latter the inner ones
\begin{multline*}
   f^{(2p+1)}_{n,n} = \frac{t^{n(p+1/2)-p(p+1)}}{p!(p+1)!\pi^{2p+1}}\int^1_0dv_1 \cdots \int^1_0dv_p
   \prod^{p}_{j=1}v^{n-2p-1}_j\left[(1-tv_j)(v^{-1}_j-1)\right]^{1/2} \Delta^2_p(v)
   \\
   \times \int^1_0du_1 \cdots \int^1_0du_{p+1} \prod^{p+1}_{j=1}u^{n-1}_j\left[(1-tu_j)(u^{-1}_j-1)\right]^{-1/2}
          \Delta^2_{p+1}(u) \prod^{p+1}_{j=1}\prod^{p}_{k=1}(u_j-t^{-1}v^{-1}_k)^{-2} .
\end{multline*}
The inner $u$-integral can be written as a determinant using the Heine identity to give
\begin{equation*}
   (p+1)! \det\left[ \int^1_0 du\, u^{n+j+k-3}\frac{[(1-tu)(u^{-1}-1)]^{-1/2}}{\prod^p_{l=1}(u-t^{-1}v^{-1}_l)^2} \right]_{1\leq j,k\leq p+1} .  
\end{equation*}
We next again employ the general identity (\ref{hankel_ID}) but because there are only $ p $ factors
in the denominator of the integrand instead of $ p+1 $ we introduce an extra factor $ (u-t^{-1}v^{-1}_{p+1})^2 $
into the numerator, which yields
\begin{multline*}
   f^{(2p+1)}_{n,n} = \frac{t^{n(p+1/2)}}{p!\pi^{2p+1}} \int^1_0dv_1 \cdots \int^1_0dv_{p}
   \prod^{p}_{j=1}v^{n-1}_j(1-tv_j)^{1/2}(v^{-1}_j-1)^{1/2} v^{2p}_{p+1}\prod^{p}_{j=1}(v_{p+1}-v_j)^{-2}
   \\ \times
   \det \left[ \begin{matrix} I(v_j,v_k) & I(v_j,v_{p+1}) \cr I(v_{p+1},v_k) & I(v_{p+1},v_{p+1}) \end{matrix} \right]_{1\leq j,k\leq p} ,
\end{multline*}
and 
\begin{equation*}
   I(v_{j},v_{k}) := \int^1_0 dx\, x^{n-1}(1-tx)^{-1/2}(x^{-1}-1)^{-1/2}\frac{(x-t^{-1}v^{-1}_{p+1})^2}{(x-t^{-1}v^{-1}_{j})(x-t^{-1}v^{-1}_{k})} .
\end{equation*}
Naturally the above integral is independent of $ v_{p+1} $ as one can verify by elementary row 
and column subtractions. Also it clearly has the structure of a coefficient in the Neumann expansion
of a Fredholm minor which is not manifestly apparent in the final result (\ref{highT_CFredholm}).
\end{proof}

\begin{remark}
The Fredholm determinant forms found here are implicit in the work of Jimbo and Miwa \cite{JM_1980}
although many details are lacking in this work for us to make a systematic comparison. Suffice it to say
that we have evaluated the Neumann expansion coefficients of their resolvent kernels, Eqs. (13) and (14), 
and found that they are essentially the same as the ones that can be deduced from (\ref{lowT_CFredholmDet}) 
and (\ref{lowT_Ckernel}).
\end{remark}

\begin{remark}
The Appell functions that appear in our application are the special cases 
\begin{equation*}
  F_{1}(n+\tfrac{1}{2};\tfrac{1}{2},1;n+1;t,tx),\qquad F_{1}(n+\tfrac{1}{2};-\tfrac{1}{2},1;n+2;t,tx) ,
\end{equation*}
of $ F_{1}(\alpha;\beta,\beta';\gamma;x,y) $
by virtue of relations amongst the parameters $ \beta+\gamma=\alpha+1 $, $ \beta'=1 $,
however we are not aware of any evidence in the literature \cite{EMOT_I}, \cite{Vi_2009}, \cite{Vi_2010} that 
they reduce to a single term of univariate hypergeometric functions or their products.
\end{remark}

\begin{remark}
Both expansions, (\ref{ffexp}) with (\ref{ff_Ieven}) and (\ref{ffhigh}) with (\ref{ff_Iodd}), are valid
for $ |t|<1 $ and are clearly inapplicable at the critical point $ t=1 $ as the multiple-integral for each individual form
factor $ f^{(.)}_{n,n} $ diverges as $ t\to 1^{-} $. The Appell function appearing in (\ref{lowT_Ckernel}) is
well-defined for $ t=1 $ ($ \gamma-\alpha-\beta=2 $) however the integrands in (\ref{NeumannLowT}) are too singular at one endpoint
for the integrals to remain finite. In contrast the Appell and hypergeometric functions appearing in
(\ref{highT_Ckernel:a}-\ref{highT_Ckernel:c}) diverge logarithmically ($ \gamma-\alpha-\beta=0 $) as
$ t\to 1$.
\end{remark}

\section{Borodin-Okounkov Identity}
\label{DFD}
\setcounter{equation}{0}
For Toeplitz determinants with symbol $ a(\theta) $ such that the Szeg\H{o} formula (\ref{szego})
holds with $ c_0=0 $, there is a general transformation identity to a Fredholm determinant
derived by Borodin and Okounkov \cite{BO_2000a} in the context of determinantal point processes.
Later it was realised that the identity had been derived previously within Toeplitz determinant
theory by Geronimo and Case \cite{GC_1979}. In \cite{Bo_2003} some special cases of the 
general identity were worked out. One of these, appropriately further specialised, is the 
Toeplitz determinant with symbol (\ref{diagonal:a}). Thus from \cite{Bo_2003} it follows in the
low temperature regime $ n\geq 0 $
\begin{equation}
   \langle \sigma_{0,0}\sigma_{n,n} \rangle = (1-t)^{1/4}\det(\mathbb{I}-\mathbb{K}_{n,n+1,\ldots}) ,
\label{KV}
\end{equation}
where $ \mathbb{K}_{n,n+1,\ldots} $ is the integral operator with kernel
\begin{multline} 
   K(i,j) = \frac{(-\tfrac{1}{2})_{i+1}(\tfrac{1}{2})_{j+1}}{i!j!}\frac{t^{(i+j)/2+1}}{(1-t)}\frac{1}{i-j}
  \\
  \times \left[ \frac{1}{j+1}{}_2F_1(-\tfrac{1}{2},\tfrac{1}{2};i+1;\frac{t}{t-1}){}_2F_1(\tfrac{1}{2},\tfrac{3}{2};j+2;\frac{t}{t-1})
               -\frac{1}{i+1}{}_2F_1(-\tfrac{1}{2},\tfrac{1}{2};j+1;\frac{t}{t-1}){}_2F_1(\tfrac{1}{2},\tfrac{3}{2};i+2;\frac{t}{t-1}) \right] ,
\label{BOkernel}
\end{multline}
supported on the successive integers $ n,n+1,n+2,\ldots $. Here $ {}_2F_1 $ is the standard 
Gauss hypergeometric function and $ K(i,i) $ is the limit $ j \to i $ of (\ref{BOkernel}).
The Kummer relation
\begin{equation*}
   {}_2F_1(a,b;c;t) = (1-t)^{-a}{}_2F_1(a,c-b;c;\frac{t}{t-1}) ,
\end{equation*}
allows (\ref{BOkernel}) to be equivalently written
\begin{multline} 
   K(i,j) = \frac{(-\tfrac{1}{2})_{i+1}(\tfrac{1}{2})_{j+1}}{i!j!}t^{(i+j)/2+1}\frac{1}{i-j}
  \\
  \times \left[ \frac{1}{j+1}{}_2F_1(\tfrac{1}{2},i+\tfrac{3}{2};i+1;t){}_2F_1(\tfrac{1}{2},j+\tfrac{1}{2};j+2;t)
               -\frac{1}{i+1}{}_2F_1(\tfrac{1}{2},j+\tfrac{3}{2};j+1;t){}_2F_1(\tfrac{1}{2},i+\tfrac{1}{2};i+2;t) \right] ,
\label{Dkernel}
\end{multline}
for $ |t|<1 $. We note that this kernel diverges as $ t\to 1 $ because of the fact that one of the Gauss
hypergeometric functions diverges with the parameter combination $ c-a-b=-1 $.

All Fredholm determinants enjoy natural $ \lambda $ generalisations. Thus according to
the general Fredholm theory \cite{WW_1958} we can expand
\begin{equation}
   \det(\mathbb{I}-\lambda^2\mathbb{K}_{n,n+1,\ldots}) 
  \\
   = 1+\sum^{\infty}_{p=1} (-\lambda^2)^p \sum_{n_1>n_2>\cdots>n_p\geq n}\det[K(n_j,n_k)]_{1\leq j,k\leq p} .
\label{KV1}
\end{equation}
This and (\ref{KV}) compared to (\ref{ffexp}) and (\ref{ffexp2}) suggest an equality between
the coefficients of $ \lambda^{2p} $ in the expansions (\ref{ffexp2}) and (\ref{KV1}).
\begin{proposition}\label{Equality} 
Let $ f^{(2p)}_{n,n} $ be specified by (\ref{ff_Ieven}) and let $ K(i,j) $ be given by (\ref{Dkernel}).
We have
\begin{equation}
   f^{(2p)}_{n,n} = (-1)^p \sum_{n_1>n_2>\cdots>n_p\geq n}\det[K(n_j,n_k)]_{1\leq j,k\leq p} .
\label{Ddet}
\end{equation}
\end{proposition} 
We will develop the proof of this statement, which is the main result of our paper, in a series of steps
starting in the next section.

A simple piece of supporting evidence in favour of (\ref{Ddet}) is that for $ p=1 $ the
leading $ t\to 0 $ behaviour of the right-hand side agrees with that implied by 
(\ref{ffsmallt}) and (\ref{ffzero:a}). This is seen by noting from (\ref{Dkernel}) that
\begin{equation*}
    K(i,i) \mathop{\sim}\limits_{t \to 0} \frac{(\tfrac{1}{2})_{i+1}(-\tfrac{1}{2})_{i+1}}{((i+1)!)^2}t^{i+1} ,
\end{equation*}
and noting furthermore that
\begin{equation*}
  \sum_{n_1\geq n}K(n_1,n_1) \mathop{\sim}\limits_{t \to 0} K(n,n) .
\end{equation*}
However, it is not straightforward to extend this analysis to show that (\ref{Ddet})
is consistent with (\ref{ffsmallt}) and (\ref{ffzero:a}) in the $ t\to 0 $ limit for $ p>1 $.

It is known from \cite{JM_1980}, \cite{JM_1980errata}, \cite{FW_2004b} that
\begin{equation}
   \sigma_{n}(t) = \left\{
   \begin{array}{lr}
     t(t-1)\frac{\displaystyle d}{\displaystyle dt}\log\langle \sigma_{0,0}\sigma_{n,n} \rangle-\frac{1}{4}t,  & T<T_c \\
     &\\
     t(t-1)\frac{\displaystyle d}{\displaystyle dt}\log\langle \sigma_{0,0}\sigma_{n,n} \rangle-\frac{1}{4}, & T>T_c
   \end{array}     \right. ,
\label{Sdefn}
\end{equation}
satisfies a differential equation which is particular case of the $\sigma$-form for Painlev\'e VI
\begin{equation}
   \left( t(t-1)\frac{d^2\sigma}{dt^2} \right)^2 = n^2\left( (t-1)\frac{d\sigma}{dt}-\sigma \right)^2
     -4\frac{d\sigma}{dt}\left( (t-1)\frac{d\sigma}{dt}-\sigma-\tfrac{1}{4} \right)\left( t\frac{d\sigma}{dt}-\sigma \right) ,
\label{sigmaform}
\end{equation}
subject to the boundary condition
\begin{equation*}
   C_{-}(n,n;\lambda) = (1-t)^{1/4} + \lambda^2\frac{(1/2)_{n}(3/2)_{n}}{4[(n+1)!]^2}t^{n+1}(1+{\rm O}(t)) ,
\end{equation*}
as $ t \to 0 $.
In \cite{BHMMOZ_2007}, \cite{MG_2010} evidence is given that the
$ \lambda $-generalisation (\ref{ffexp2}) satisfies the same differential equation, with the
$ \lambda $-dependence entering only in the boundary condition.

On the other hand the identity (\ref{KV}) tells us that 
$ \log\det(\mathbb{I}-\mathbb{K}_{n,n+1,\ldots}) $ satisfies a simple variant of (\ref{sigmaform}).
Proposition \ref{Equality} together with the validity of the conjecture from \cite{MG_2010} then 
gives that $ \log\det(\mathbb{I}-\lambda^2\mathbb{K}_{n,n+1,\ldots}) $
satisfies the same differential equation up to the boundary condition. Specifically, this 
would imply that with
\begin{equation}
   \sigma(t) = t(t-1)\frac{d}{dt} \log\det(\mathbb{I}-\lambda^2\mathbb{K}_{n,n+1,\ldots}) ,
\label{SigmaDefn}
\end{equation}
the differential equation (\ref{sigmaform}) is satisfied and furthermore that this equation
together with the boundary condition
\begin{equation*}
   \sigma(t) \mathop{\sim}\limits_{t\to 0} -\lambda^2 (n+1)\left. c_0\right|_{p=1}t^{n+1} ,
\end{equation*}
where $ c_0 $ is specified by (\ref{ffzero:a}), completely determines
$ \det(\mathbb{I}-\lambda^2\mathbb{K}_{n,n+1,\ldots}) $. It remains an open problem to verify this
characterisation.

\section{Bi-orthogonal Polynomials on the Unit Circle}
\label{BOPS}
\setcounter{equation}{0}

Our first step in the derivation of (\ref{Ddet}) is to build up a theory of bi-orthogonal polynomials
on the unit circle and its relationship to general, non-hermitian (i.e. complex weight) Toeplitz
determinants. The study by Geronimo and Case \cite{GC_1979}
analysed systems of orthogonal polynomials on the unit circle, or equivalently the case of hermitian Toeplitz 
determinants and is therefore inadequate for our purposes.
However we will see that much of their analysis generalises in an obvious manner provided one makes suitable 
distinctions between variables, most notably splitting the orthogonal system into a bi-orthogonal one, 
and remains valid by replacing certain positivity requirements by non-vanishing ones.
In his 1966 prescient study of Toeplitz determinants applied to the spin correlations of the Ising
model Wu \cite{Wu_1966} uncovered some of the structures that exist in systems of bi-orthogonal
polynomials without explicitly appreciating or exploiting that fact. Much of the theoretical 
development outlined in this section of our present study draws upon the foundations given in \cite{FW_2006a}.
 
We consider a complex function for our weight $ w(z) $, analytic in the cut
complex $ z $-plane and which possesses a Fourier expansion
\begin{equation*}
  w(z) = \sum_{k=-\infty}^{\infty} w_{k}z^k, \quad
  w_{k} = \int_{\T} \frac{d\zeta}{2\pi i\zeta} w(\zeta)\zeta^{-k},
\end{equation*}
where $ z \in D \subset \C $ and $ \T $ denotes the unit circle $ |\zeta|=1 $ with 
$ \zeta=e^{i\theta}, \theta \in (-\pi,\pi] $. 
Hereafter we will assume that $ z^jw(z), z^jw'(z) \in L^1(\T) $ for all 
$ j \in \Z $. We will also assume that the trigonometric sum converges in an annulus
$ D=\{z\in \C:\Delta_1 < |z| < \Delta_2 \} $ and $ \T \subset D $. The doubly infinite 
sequence $ \{ w_k \}^{\infty}_{k=-\infty} $ are the trigonometric moments of the 
distribution $ w(e^{i\theta})d\theta/2\pi $ and define the trigonometric moment 
problem (see \cite{Simon_I_2005} for an up-to-date account of the hermitian case $ w^{\dagger}_p=w_{-p} $). 
In \cite{GC_1979} their variable $ Z = z=e^{i\theta} $ and their measure
$ d\rho(\theta) $ is related to our weight by
\begin{equation*}
   d\rho(\theta) = w(z)\frac{dz}{iz} ,
\end{equation*}
and consequently their moments $ c_{n} = w_{n} $. In the work of Wu \cite{Wu_1966}
his weight $ C(z) = w(z) $ and so his moments $ c_{n} = w_{n} $. In the Borodin-Okounkov work
we identify their weight $ \phi(z)=w(-z) $ and $ V(z)=\log w(-z) $.
Utilising the trigonometric moments we define the Toeplitz determinants $ n\geq 1 $
\begin{align}
   I_{n}[\zeta^{\epsilon}w(\zeta)]
  & \mathrel{\mathop:}= \det \left[ \int_{\T} \frac{d\zeta}{2\pi i\zeta} w(\zeta)\zeta^{\epsilon-j+k}
           \right]_{0 \leq j,k \leq n-1}
     = \det \left[ w_{-\epsilon+j-k} \right]_{0 \leq j,k \leq n-1},
\label{Tdet} \\
  & = \frac{1}{n!} \int_{\T^n}
      \prod^{n}_{l=1} \frac{d\zeta_l}{2\pi i\zeta_l} w(\zeta_l)\zeta_l^{\epsilon}
      \prod_{1 \leq j<k \leq n} |\zeta_{j}-\zeta_{k}|^2,
\end{align}
where $ \epsilon $ will take the integer values $ 0,\pm 1 $.
In \cite{GC_1979} their Toeplitz determinant $ D_{n} = I_{n+1} $.

We define a system of bi-orthogonal polynomials 
$ \{ \varphi_n(z),\bar{\varphi}_n(z) \}^{\infty}_{n=0} $ with respect to the 
weight $ w(z) $ on the unit circle by the orthogonality relation
\begin{equation}
  \int_{\T} \frac{d\zeta}{2\pi i\zeta} w(\zeta)\varphi_m(\zeta)\bar{\varphi}_n(\bar{\zeta})  
   = \delta_{m,n} .
\label{ops_onorm}
\end{equation}
This system is taken to be orthonormal and the coefficients in a monomial basis 
are defined by
\begin{align*}
   \varphi_n(z)
   & = \kappa_n z^n + l_n z^{n-1}+ m_n z^{n-2} + \ldots + \varphi_n(0)
     = \sum^{n}_{j=0} c_{n,j}z^j,
   \\
   \bar{\varphi}_n(z)
   & = \bar{\kappa}_n z^n + \bar{l}_n z^{n-1}+ \bar{m}_n z^{n-2} + \ldots + \bar{\varphi}_n(0)
     = \sum^{n}_{j=0} \bar{c}_{n,j}z^j,
\end{align*}
where $ \bar{\kappa}_n $ is chosen to be equal to $ \kappa_n $ without loss of generality
(this has the effect of rendering many results formally identical to the pre-existing
theory of orthogonal polynomials). In \cite{GC_1979} the variable $ K(n) = \kappa_{n} $.
Notwithstanding the notation $ \bar{c}_{n,j} $ 
in general is not equal to the complex conjugate of $ c_{n,j} $ and is independent of it.
We also define the reverse or reciprocal polynomial by
\begin{equation*}
   \varphi^{*}_n(z) := z^n\bar{\varphi}_n(1/z) = \sum^{n}_{j=0} \bar{c}_{n,j}z^{n-j} .
\end{equation*}
The bi-orthogonal polynomials can equivalently be defined up to an overall factor by the orthogonality
with respect to the monomials
\begin{equation}
  \int_{\T} \frac{d\zeta}{2\pi i\zeta} w(\zeta)\varphi_n(\zeta)\overline{\zeta^j}  
   = 0 \qquad 0 \leq j \leq n-1 ,
\label{ops_orthog:a}
\end{equation}
whereas their reciprocal polynomials can be similarly defined by
\begin{equation}
  \int_{\T} \frac{d\zeta}{2\pi i\zeta} w(\zeta)\varphi^*_n(\zeta)\overline{\zeta^j}  
   = 0 \qquad 1 \leq j \leq n .
\label{ops_orthog:b}
\end{equation}
The polynomials defined in \cite{GC_1979} are related to ours by
\begin{equation*}
  \phi(Z,n) = \varphi_{n}(z), \quad \bar{\phi}(Z,n) = \bar{\varphi}_{n}(z), \quad \phi^{*}(Z,n) = \varphi^{*}_{n}(z) . 
\end{equation*}
In \cite{Wu_1966} his Eqs. (2.7), (2.12) are directly comparable with our orthogonality condition 
(\ref{ops_orthog:b}), and consequently his polynomial $ X(z) = \kappa_n\varphi^{*}_{n}(z) $ after
reconciling the normalisations.

The linear system of equations for the coefficients $ c_{n,j}, \bar{c}_{n,j} $
arising from 
\begin{align}
  \bar{c}_{n,n}\int_{\T} \frac{d\zeta}{2\pi i\zeta} w(\zeta)\varphi_n(\zeta)\bar{\zeta}^m  
 & = \begin{cases} 0 \quad 0 \leq m \leq n-1 \\ 1 \quad m=n \end{cases} , 
  \label{ops_CoeffEqn:a}\\
  c_{n,n}\int_{\T} \frac{d\zeta}{2\pi i\zeta} w(\zeta)\zeta^m\bar{\varphi}_n(\bar{\zeta})
 & = \begin{cases} 0 \quad 0 \leq m \leq n-1 \\ 1 \quad m=n \end{cases} ,
  \label{ops_CoeffEqn:b}
\end{align}
has the solution
\begin{align*}
  c_{nj} & = \frac{1}{\bar{c}_{n,n}}
         \frac{\det \begin{pmatrix}
                     w_{0} 	& \ldots & 0 	& \ldots	& w_{-n} \cr
                     \vdots	& \vdots & \vdots & \vdots 	& \vdots \cr
                     w_{n-1}	& \ldots & 0 	& \ldots 	& w_{-1} \cr
                     w_{n} 	& \ldots & 1 	& \ldots 	& w_{0}  \cr
                    \end{pmatrix} }
              {\det \begin{pmatrix}
                     w_{0} & \ldots & w_{-n} \cr
                     \vdots& \vdots & \vdots \cr
                     w_{n} & \ldots & w_{0}  \cr
                    \end{pmatrix} } ,
  \\
  \bar{c}_{nj} & = \frac{1}{c_{n,n}}
         \frac{\det \begin{pmatrix}
                     w_{0} 	& \ldots & 0 	& \ldots 	& w_{n}  \cr
                     \vdots	& \vdots & \vdots & \vdots 	& \vdots \cr
                     w_{-n+1} 	& \ldots & 0 	& \ldots 	& w_{1}  \cr
                     w_{-n} 	& \ldots & 1 	& \ldots 	& w_{0}  \cr
                    \end{pmatrix} }
              {\det \begin{pmatrix}
                     w_{0} & \ldots & w_{n} \cr
                     \vdots& \vdots & \vdots   \cr
                     w_{-n}& \ldots & w_{0} \cr
                    \end{pmatrix} } ,
\end{align*}
and in particular one has the following results.
\begin{proposition}[\cite{Ba_1960},\cite{FW_2006a}]
The leading and trailing coefficients of the polynomials $ \varphi_n(z) $, $ \bar{\varphi}_n(z) $
are 
\begin{gather*}
   c_{nn} =  \bar{c}_{nn} = \kappa_{n} 
          = \frac{1}{\kappa_n}\frac{I_{n}[w(\zeta)]}{I_{n+1}[w(\zeta)]}, \\
   c_{n0} = \varphi_{n}(0) = (-1)^n\frac{1}{\kappa_n}\frac{I_{n}[\zeta w(\zeta)]}{I_{n+1}[w(\zeta)]}, \quad
   \bar{c}_{n0} = \bar{\varphi}_{n}(0) 
          = (-1)^n\frac{1}{\kappa_n}\frac{I_{n}[\zeta^{-1} w(\zeta)]}{I_{n+1}[w(\zeta)]} .
\end{gather*}
\end{proposition}
The following existence theorem for the bi-orthogonal polynomial system is due to G. Baxter.
\begin{proposition}[\cite{Ba_1960}]
The bi-orthogonal polynomial system $ \{\varphi_{n}(z),\varphi^*_{n}(z)\}^{\infty}_{n=0} $ exists
if and only if $ I_{n}[w(\zeta)] \neq 0 $ for all $ n \geq 1 $. 
\end{proposition}

The system is alternatively defined by the sequence of ratios $ r_n = \varphi_n(0)/\kappa_n $,
known as reflection coefficients because of their role in the scattering theory
formulation of orthogonal polynomials on the unit circle, together with a companion quantity
$ \bar{r}_n = \bar{\varphi}_n(0)/\kappa_n $. As in the Szeg\H{o} theory 
\cite{ops_Sz} $ r_n $ and $ \bar{r}_n $ are related to the above Toeplitz 
determinants by
\begin{equation*}
  r_n = (-1)^n\frac{I_{n}[\zeta w(\zeta)]}{I_{n}[w(\zeta)]}, \quad
  \bar{r}_n = (-1)^n\frac{I_{n}[\zeta^{-1} w(\zeta)]}{I_{n}[w(\zeta)]}.
\end{equation*}
In \cite{GC_1979} we have the correspondences of $ \alpha(n) = \varphi_{n}(0) $,
$ \overline{\alpha(n)} = \bar{\varphi}_{n}(0) $, $ a(n) = \kappa_{n+1}/\kappa_{n} $ and $ b(n) = r_{n} $.
The Toeplitz determinants of central interest can then be recovered through the 
following result, which was given by Baxter in 1961 \cite{Ba_1961} and is the generalisation of
Eq. (II.16) of Geronimo and Case \cite{GC_1979}.
\begin{proposition}[\cite{Ba_1961},\cite{FW_2006a}]
The sequence of $ \{I_{n}[w]\}^{\infty}_{n=0} $ satisfy the recurrence
\begin{equation}
   \frac{I_{n+1}[w]I_{n-1}[w]}{(I_{n}[w])^2}
   = 1 - r_{n}\bar{r}_n, \quad n \geq 1 ,
\label{ops_I0}
\end{equation}
subject to the condition $ r_{n}\bar{r}_n \neq 1 $ for all $ n \geq 1 $. 
\end{proposition}
Further identities from the Szeg\H{o} theory that generalise are those that relate 
the leading coefficients back to the reflection coefficients. For example we have
\begin{equation}
   \kappa_n^2 = \kappa_{n-1}^2 + \varphi_n(0)\bar{\varphi}_n(0) ,
\label{ops_kappa}
\end{equation}
which is an extension of Eq. (II.10) in \cite{GC_1979}.

A fundamental consequence of the orthogonality conditions is a system of coupled linear first order
difference equations. They can constitute the starting point for the theory of a bi-orthogonal polynomial system 
rather than the orthogonality conditions and this was how Baxter \cite{Ba_1960}, \cite{Ba_1961} developed his 
theory. These recurrence relations are equivalent to Eq. (II.7, 8) of Geronimo and Case \cite{GC_1979}.
\begin{proposition}[\cite{Ba_1960,Ba_1961}]
The polynomial pair $ \varphi_{n}(z), \varphi^*_n(z) $ satisfy the coupled recurrence relations
\begin{align}
  \kappa_n  \varphi_{n+1}(z)
   & = \kappa_{n+1}z \varphi_{n}(z)+\varphi_{n+1}(0) \varphi^*_n(z) ,
  \label{ops_rr:a} \\
  \kappa_n \varphi^*_{n+1}(z)
   & = \kappa_{n+1} \varphi^*_{n}(z)+\bar{\varphi}_{n+1}(0) z\varphi _n(z) .
  \label{ops_rr:b}
\end{align}
\end{proposition}

Together equations (\ref{ops_rr:a},\ref{ops_rr:b}) are equivalent to a single second order linear difference 
equation and therefore admit another linearly independent solution. We define another polynomial solution 
pair $ \psi_n(z), \psi^*_n(z) $ - the polynomial of the second kind or associated polynomial
\begin{equation}
  \psi_n(z) 
  := \int_{\T}\frac{d\zeta}{2\pi i\zeta}\frac{\zeta+z}{\zeta-z}w(\zeta)
               [\varphi_n(\zeta)-\varphi_n(z)],
       \quad n \geq 1, \quad \psi_0 := \kappa_0w_0 = 1/\kappa_0,
\label{ops_psi:a}
\end{equation}
and its reciprocal polynomial $ \psi^*_n(z) $. The integral formula for $ \psi^*_{n} $ is
\begin{equation}
  \psi^*_n(z) 
  := -\int_{\T}\frac{d\zeta}{2\pi i\zeta}\frac{\zeta+z}{\zeta-z}w(\zeta)
                [z^n\bar{\varphi}_n(\bar{\zeta})-\varphi^*_n(z)],
       \quad n \geq 1, \quad \psi^*_0 := 1/\kappa_0.
\label{ops_psi:b}
\end{equation}
A central object in our theory is the Carath\'eodory function - the generating function of the Toeplitz elements
\begin{equation}
   F(z) := \int_{\T}\frac{d\zeta}{2\pi i\zeta}\frac{\zeta+z}{\zeta-z}w(\zeta) ,\quad z \notin \T,
\label{ops_Cfun}
\end{equation}
which defines the inner and outer functions with expansions inside and outside the unit circle
respectively
\begin{align*}
   F^{<}(z) & = w_0 + 2 \sum^{\infty}_{k=1}w_{k} z^k, \quad \text{if $ |z| < 1 $}, \\
   F^{>}(z) & = -w_0 - 2 \sum^{\infty}_{k=1}w_{-k} z^{-k}, \quad \text{if $ |z| > 1 $}.
\end{align*}
Having these definitions one requires two non-polynomial solutions 
$ \epsilon_n(z), \epsilon^*_n(z) $ to the recurrences (\ref{ops_rr:a},\ref{ops_rr:b}) and these are constructed 
as linear combinations of the polynomial solutions according to 
\begin{align}
   \epsilon_n(z) := \psi_n(z)+F(z)\varphi_n(z)
   & = \int_{\T}\frac{d\zeta}{2\pi i\zeta}\frac{\zeta+z}{\zeta-z}w(\zeta)
                   \varphi_n(\zeta) ,
   \label{ops_eps:a} \\
   \epsilon^*_n(z) := \psi^*_n(z)-F(z)\varphi^*_n(z)
   & = \frac{1}{\kappa_n} 
          -\int_{\T}\frac{d\zeta}{2\pi i\zeta}\frac{\zeta+z}{\zeta-z}w(\zeta)
                   \varphi^*_n(\zeta) ,
   \label{ops_eps:b}
\end{align}
for $ n\geq 1 $ and $ \epsilon_0 = \kappa_0(w_0+F) $, $ \epsilon^*_0 = \kappa_0(w_0-F) $.
Equation (\ref{ops_eps:b}) along with the definition (\ref{ops_orthog:b}), is directly comparable to
Eqs. (2.10), (2.14) and to Eqs. (2.11), (2.15) of the Wu study \cite{Wu_1966} in their
respective $z$-regimes.
\begin{proposition}[\cite{Ge_1961},\cite{Ge_1962},\cite{Ge_1977},\cite{JNT_1989},\cite{FW_2006a}]
The associated functions $ \epsilon_{n}(z), \epsilon^*_{n}(z) $ satisfy a variant of 
the coupled recurrence relations (\ref{ops_rr:a},\ref{ops_rr:b}) namely
\begin{align}
  \kappa_n  \epsilon_{n+1}(z)
   & = \kappa_{n+1}z \epsilon_{n}(z)-\varphi_{n+1}(0) \epsilon^*_n(z) ,
  \label{ops_rre:a} \\
  \kappa_n \epsilon^*_{n+1}(z)
   & = \kappa_{n+1} \epsilon^*_{n}(z)-\bar{\varphi}_{n+1}(0) z\epsilon _n(z) .
  \label{ops_rre:b}
\end{align}
\end{proposition}
The linear independence of the two solution sets to the coupled recurrences (\ref{ops_rr:a},\ref{ops_rr:b})
has the following consequences.
\begin{proposition}[\cite{Ge_1961},\cite{FW_2006a}]
The Casoratians of the polynomial solutions $ \varphi_n, \varphi^*_n, \psi_n, \psi^*_{n} $ or of
the polynomial and non-polynomial solutions $ \varphi_n, \varphi^*_n, \epsilon_n, \epsilon^*_{n} $ are
\begin{align}
   \varphi_{n+1}(z)\psi_n(z) - \psi_{n+1}(z)\varphi_n(z) 
  & = \varphi_{n+1}(z)\epsilon_n(z) - \epsilon_{n+1}(z)\varphi_n(z)
    = 2\frac{\varphi_{n+1}(0)}{\kappa_n}z^n ,
  \label{ops_Cas:a} \\
   \varphi^*_{n+1}(z)\psi^*_n(z) - \psi^*_{n+1}(z)\varphi^*_n(z) 
  & = \varphi^*_{n+1}(z)\epsilon^*_n(z) - \epsilon^*_{n+1}(z)\varphi^*_n(z) 
    = 2\frac{\bar{\varphi}_{n+1}(0)}{\kappa_n}z^{n+1} , 
  \label{ops_Cas:b} \\
   \varphi_{n}(z)\psi^*_n(z) + \psi_{n}(z)\varphi^*_n(z) 
  & = \varphi_{n}(z)\epsilon^*_n(z) + \epsilon_{n}(z)\varphi^*_n(z)
    = 2z^n .  
  \label{ops_Cas:c}
\end{align} 
\end{proposition}

We note that the recurrence relations for the associated functions 
$ \epsilon_n(z) $, $ \epsilon^*_n(z) $ given in (\ref{ops_rre:a},\ref{ops_rre:b}) differ
from those of the polynomial systems (\ref{ops_rr:a},\ref{ops_rr:b}) by a reversal of
the signs of $ \varphi_n(0) $, $ \bar{\varphi}_n(0) $. We can compensate for this by constructing
the $ 2\times 1 $ vectors
\begin{equation}
   Y_n(z) := 
   \begin{pmatrix}
          \varphi_n(z)   \cr
          \varphi^*_n(z)
   \end{pmatrix} \quad\text{or}\quad
   \begin{pmatrix}
           \epsilon_n(z) \cr
          -\epsilon^{\vphantom{I}*}_n(z)
   \end{pmatrix} ,
\label{ops_Ydefn}
\end{equation}
which can be directly compared with the polynomial solution $ \Psi^*(Z,n) = Y_n(z) $ in \cite{GC_1979}.

Our extension of the matrix recurrence relations Eq. (II.12, 13) of Geronimo and Case \cite{GC_1979}
is the following result.
\begin{corollary}[\cite{FW_2006a}]
The recurrence relations for a general system of bi-orthogonal polynomials 
(\ref{ops_rr:a},\ref{ops_rr:b}) and their associated functions 
(\ref{ops_rre:a},\ref{ops_rre:b}) are equivalent to the first order matrix recurrence
\begin{equation}
   Y_{n+1} = K_n Y_{n}
  := \frac{1}{\kappa_n}
       \begin{pmatrix}
              \kappa_{n+1} z   & \varphi_{n+1}(0) \cr
              \bar{\varphi}_{n+1}(0) z & \kappa_{n+1} \cr
       \end{pmatrix} Y_{n} . 
\label{ops_Yrecur}
\end{equation} 
A consequence of (\ref{ops_kappa}) is that $ K_n $ has the property $ \det K_n = z $.
\end{corollary}
In Geronimo and Case \cite{GC_1979} the matrix $ A^{*}(n) = K_{n} $. In addition to this matrix they 
defined the matrix
\begin{equation*}
  A(n) = \begin{pmatrix}
    1 & 0 \cr
    0 & z^{-n-1}
  \end{pmatrix}
  K_n
  \begin{pmatrix}
    z^{n} & 0 \cr
    0 & 1
  \end{pmatrix} ,
\end{equation*}
which figures in an alternative matrix recurrence of the matrix variable 
\begin{equation*}
  \Psi(Z,n) = \begin{pmatrix}
                    \varphi_{n}(z) \cr
                    z^{-n}\varphi^{*}_{n}(z) 
              \end{pmatrix} .
\end{equation*}
Accordingly we can also define an alternative matrix of polynomial solutions $ \Psi_{n} $ as
\begin{equation}
   \Psi_{n}(z) = \begin{pmatrix} \varphi_{n}(z) \cr \bar{\varphi}_{n}(z^{-1}) \end{pmatrix} ,
\label{PsiDefn}
\end{equation}
which satisfies a variant of (\ref{ops_Yrecur}), namely 
\begin{equation}
   \Psi_{n+1}
   = \frac{1}{\kappa_n}
       \begin{pmatrix}
              \kappa_{n+1} z   & \varphi_{n+1}(0)z^{n} \cr
              \bar{\varphi}_{n+1}(0)z^{-n} & \kappa_{n+1}z^{-1} \cr
       \end{pmatrix}\Psi_{n} .
\label{PsiRecur}
\end{equation} 
The utility of this definition is that the determinant of the matrix given above is unity.

The analogue of the Christoffel-Darboux summation formula is given by the following result.
\begin{proposition}[\cite{Ba_1961},\cite{FW_2006a}]
The summation identity
\begin{equation}
  \sum^{n}_{j=0} \varphi_j(z)\bar{\varphi}_j(\bar{\zeta}) 
  =
  \frac{\varphi^*_n(z)\overline{\varphi^*_n}(\bar{\zeta})
        -z\bar{\zeta}\varphi_n(z)\bar{\varphi}_n(\bar{\zeta})}{1-z\bar{\zeta}}
  = 
  \frac{\varphi^*_{n+1}(z)\overline{\varphi^*_{n+1}}(\bar{\zeta})
        -\varphi_{n+1}(z)\bar{\varphi}_{n+1}(\bar{\zeta})}{1-z\bar{\zeta}} ,
\label{ops_CD}
\end{equation}
holds for $ z\bar{\zeta} \not= 1 $ and $ n \geq 0 $. Here 
\begin{equation*}
   \overline{\varphi^*_n}(\bar{\zeta}) = \bar{\zeta}^n\varphi_n(1/\bar{\zeta}) .
\end{equation*}
\end{proposition}
Another related bilinear identity is the formula for the discrete Wronskian.
\begin{proposition}[\cite{GC_1979}]\label{wronskian}
The discrete Wronskian of two solutions $ \Psi_{1,n}, \Psi_{2,n} $ to the recurrence 
system (\ref{PsiRecur}) is
\begin{equation*}
    W[\Psi_{1,n}, \Psi_{2,n}] = \Psi_{1,n}^{T}
    \begin{pmatrix} 0 & -1 \cr 1 & 0 \end{pmatrix}
                                \Psi_{2,n} ,
\end{equation*}
i.e. $ W[\Psi_{1,n+1}, \Psi_{2,n+1}] = W[\Psi_{1,n}, \Psi_{2,n}] $. Naturally this latter relation still holds
for non-polynomial solutions of (\ref{PsiRecur}) as well.
\end{proposition}

Because we will be expanding the recurrence relation solutions about $ z=0, \infty $
we give the leading order terms in the expansions of 
$ \varphi_n(z), \varphi^*_n(z), \epsilon_n(z), \epsilon^*_{n}(z) $ both inside and 
outside the unit circle.
\begin{corollary}[\cite{FW_2006a}]
The bi-orthogonal polynomials $ \varphi_n(z), \varphi^*_n(z) $ have the expansions
\begin{align}
   \varphi_n(z) & = 
   \begin{cases}
      \varphi_n(0)
       + \dfrac{1}{\kappa_{n-1}}
         (\kappa_n\varphi_{n-1}(0)+\varphi_n(0)\bar{l}_{n-1})z
       + {\rm O}(z^2) & |z| < 1 ,\\
      \kappa_n z^n + l_n z^{n-1} + {\rm O}(z^{n-2}) & |z| > 1 ,
   \end{cases} \label{ops_phiexp:a} \\
   \varphi^*_n(z) & = 
   \begin{cases}
      \kappa_n + \bar{l}_n z + {\rm O}(z^{2}) & |z| < 1 ,\\
      \bar{\varphi}_n(0) z^n
       + \dfrac{1}{\kappa_{n-1}}
         (\kappa_n\bar{\varphi}_{n-1}(0)+\bar{\varphi}_n(0)l_{n-1})z^{n-1}
       + {\rm O}(z^{n-2}) & |z| > 1 .\\
   \end{cases} \label{ops_phiexp:b}
\end{align}
The associated functions have the expansions
\begin{equation}
   \dfrac{\kappa_n}{2}\epsilon^{<}_n(z) = 
      z^n - \dfrac{\bar{l}_{n+1}}{\kappa_{n+1}}z^{n+1} + {\rm O}(z^{n+2}), \quad |z| < 1 ,
\label{ops_epsexp:a}
\end{equation}
\begin{equation}
   \dfrac{\kappa_n}{2}\epsilon^{>}_n(z) = 
         \dfrac{\varphi_{n+1}(0)}{\kappa_{n+1}}z^{-1}
       + \left(\dfrac{\kappa^2_n}{\kappa^2_{n+1}} \dfrac{\varphi_{n+2}(0)}{\kappa_{n+2}} -\dfrac{\varphi_{n+1}(0)}{\kappa_{n+1}} \dfrac{l_{n+1}}{\kappa_{n+1}} \right)z^{-2}
       + {\rm O}(z^{-3}), \quad |z| > 1 ,
\label{ops_epsexp:b}
\end{equation}
\begin{equation}
   \dfrac{\kappa_n}{2}\epsilon^{*<}_n(z) = 
         \dfrac{\bar{\varphi}_{n+1}(0)}{\kappa_{n+1}}z^{n+1} + \left(\dfrac{\kappa^2_n}{\kappa^2_{n+1}} \dfrac{\bar{\varphi}_{n+2}(0)}{\kappa_{n+2}}
               -\dfrac{\bar{\varphi}_{n+1}(0)}{\kappa_{n+1}} \dfrac{\bar{l}_{n+1}}{\kappa_{n+1}} \right)z^{n+2}
       + {\rm O}(z^{n+3}), \quad |z| < 1 ,
\label{ops_epsexp:c}
\end{equation}
\begin{equation}
   \dfrac{\kappa_n}{2}\epsilon^{*>}_n(z) = 
       1 - \dfrac{l_{n+1}}{\kappa_{n+1}}z^{-1} + \left( \dfrac{l_{n+2}l_{n+1}}{\kappa_{n+2}\kappa_{n+1}} -\dfrac{m_{n+2}}{\kappa_{n+2}}
         \right)z^{-2} + {\rm O}(z^{-3}), \quad |z| > 1 .
\label{ops_epsexp:d}
\end{equation}
\end{corollary}
In the notations of Wu \cite{Wu_1966} we make the identifications
\begin{equation*}
  2U(z) = -\delta_{n,0}-\kappa_n z^{-n}\epsilon^{*<}_{n}(z), \qquad
  2V(z^{-1}) = -2+\kappa_n \epsilon^{*>}_{n}(z) .
\end{equation*}

Two crucial relations are the jump conditions on the associated functions at $ |z|=1 $.
\begin{proposition}\label{RHjump}
The associated functions defined in the interior of the unit circle $ \epsilon^{<}_{n}(z), \epsilon^{*<}_{n}(z) $
differ from those in the exterior $ \epsilon^{>}_{n}(z), \epsilon^{*>}_{n}(z) $ on $ |z|=1 $ through
the jump conditions
\begin{gather}
   w(z)\varphi_{n}(z) = -\tfrac{1}{2}\epsilon^{>}_{n}(z)+\tfrac{1}{2}\epsilon^{<}_{n}(z) ,
\label{jump} \\
   w(z)\varphi^{*}_{n}(z) = \tfrac{1}{2}\epsilon^{*>}_{n}(z)-\tfrac{1}{2}\epsilon^{*<}_{n}(z) .
\label{jump*}
\end{gather}
\end{proposition}
The latter relation (\ref{jump*}) is equivalent to Eq. (2.16) of \cite{Wu_1966} which states that
\begin{equation*} 
   C(\zeta)X(\zeta) = 1+\zeta^n U(\zeta)+V(\zeta^{-1}) ,
\end{equation*}
on $ |\zeta|=1 $.

We make the observation on the weight (\ref{diagonal:a}) with $ T>T_c $ or $ \K<1 $ that
\begin{equation*}
   w^{T>T_c}(z) = \left( \frac{1-\K^{-1}z^{-1}}{1-\K^{-1}z} \right)^{1/2} = \left( \frac{1}{z^2}\frac{\K z-1}{\K z^{-1}-1} \right)^{1/2}
                = -z^{-1}\left( \frac{1-\K z}{1-\K z^{-1}} \right)^{1/2} = -z^{-1}\tilde{w}(z) .
\end{equation*}
Clearly $ w^{T>T_c}(z) $ has non-zero winding number and so Szeg\H{o}'s theorem is not
applicable, however $ \tilde{w} $ has zero winding number and is given by the corresponding weight for $ T<T_c $ with
$ \K \mapsto \K^{-1} $ and $ z \mapsto z^{-1} $. This latter transformation simply implies
$ w_{n} \leftrightarrow w_{-n} $
and the Toeplitz determinant is unaffected. However the modification of the weight by a rational
factor $ w(z) \mapsto z^{-1}w(z) $ induces a Christoffel-Uvarov-Geronimus transformation.
The elementary Christoffel-Uvarov-Geronimus transformation and it inverse 
$ w \mapsto z^{\pm 1}w $ are known from \cite{IR_1992}, \cite{Wi_2009b} - for $ w \mapsto zw $ this
is the $ K=1, L=0, \alpha\to 0 $ case in \cite{Wi_2009b}
\begin{gather*}
  I_{n} \mapsto I^{+}_{n} = (-1)^nr_{n}I_{n} ,
  \\
  \kappa_n \mapsto \kappa^{+}_{n}, \quad (\kappa^{+}_{n})^2 = -\kappa^2_{n}\frac{r_{n}}{r_{n+1}} ,
  \\
  r_{n} \mapsto r^{+}_{n} = r_{n}-\frac{\kappa^2_{n-1}}{\kappa^2_{n}}\frac{r_{n+1}r_{n-1}}{r_{n}} ,
  \\ 
  \bar{r}_{n} \mapsto \bar{r}^{+}_{n} = \frac{1}{r_{n}} ,
\end{gather*}
whilst for $ w \mapsto z^{-1}w $ this is the $ K=0, L=1, \beta\to 0 $ case
\begin{gather*}
  I_{n} \mapsto I^{-}_{n} = (-1)^n\bar{r}_{n}I_{n} ,
  \\
  \kappa_n \mapsto \kappa^{-}_{n}, \quad (\kappa^{-}_{n})^2 = -\kappa^2_{n}\frac{\bar{r}_{n}}{\bar{r}_{n+1}} ,
  \\
  r_{n} \mapsto r^{-}_{n} = \frac{1}{\bar{r}_{n}} ,
  \\ 
  \bar{r}_{n} \mapsto \bar{r}^{-}_{n} = \bar{r}_{n}-\frac{\kappa^2_{n-1}}{\kappa^2_{n}}\frac{\bar{r}_{n+1}\bar{r}_{n-1}}{\bar{r}_{n}} .
\end{gather*}

\section{Scattering States and Jost Functions}
\label{ST}
\setcounter{equation}{0}

We will make a number of assumptions, some of which we place on the symbol and some on reflection coefficients 
directly, in order to proceed with our analysis -\\
\noindent(i)
The symbol $ w(z) $ is absolutely integrable in the sense $ w(e^{i\theta}) \in L_p[-\pi,\pi] $ for $ p>1 $,\\
\noindent(ii)
That the symbol has zero winding number and $ \log w(z) $ satisfies the condition of Szeg\H{o}'s theorem
\begin{equation*}
  \sum^{\infty}_{p=-\infty}|p||c_{p}|^2 < \infty ,
\end{equation*}
\noindent(iii)
Corresponding to Eq. (III.18) of \cite{GC_1979} we assume 
that
\begin{equation*}
   \sum^{\infty}_{n=0} \sqrt{|r_{n}\bar{r}_{n}|} < \infty .
\end{equation*}
\noindent(iv)
In addition we assume $ \lim_{n \to \infty}\kappa_{n} = \kappa_{\infty} $ so that
\begin{equation*}
   \lim_{n \to \infty}\frac{\kappa_{n+1}}{\kappa_{n}} = 1 ,
\end{equation*}
and as a consequence we have
\begin{equation*}
   \lim_{n \to \infty}A(n) = 
   \begin{pmatrix}
    z & 0 \\
    0 & z^{-1}
   \end{pmatrix}
   = A(\infty) .
\end{equation*}
One can check that these facts are verified in our application to the Ising diagonal correlations in
an a posteriori demonstration, however we do not pursue those details here. 

\begin{definition}
Let $ h(z) $ have a Laurent expansion in $ z $. We define the non-negative, positive and negative parts 
of $ h $ as
\begin{equation*}
    \left[ h \right]_{\geq 0} \mathrel{\mathop:}= \sum^{\infty}_{n=0} h_n z^n, \quad
    \left[ h \right]_{> 0} \mathrel{\mathop:}= \sum^{\infty}_{n=1} h_n z^n, \quad
    \left[ h \right]_{< 0} \mathrel{\mathop:}= \sum^{-1}_{n=-\infty} h_n z^n .
\end{equation*}
\end{definition}

\begin{definition}[\cite{GC_1979}]\label{SSdefn}
Let us define the scattering states for $ n\geq 0 $ by
\begin{equation}
    \Psi_{+n}(z) \mathrel{\mathop:}= \begin{pmatrix} \phi_{+n}(z) \cr \hat{\phi}_{+n}(z) \end{pmatrix}, \quad
    \Psi_{-n}(z) \mathrel{\mathop:}= \begin{pmatrix} \hat{\phi}_{-n}(z) \cr \phi_{-n}(z) \end{pmatrix} .
\label{scatterS}
\end{equation}
They satisfy (\ref{PsiRecur}) subject to the boundary conditions
\begin{align}
   \lim_{n\to \infty} |\phi_{+n}-z^{n}| & = 0, \quad |z|<1 ,
\label{ssBC:a}   \\
   \lim_{n\to \infty} |\hat{\phi}_{+n}| & = 0, \quad |z|<1 ,
\label{ssBC:b}   \\
   \lim_{n\to \infty} |\hat{\phi}_{-n}| & = 0, \quad |z|>1 ,
\label{ssBC:c}   \\
   \lim_{n\to \infty} |\phi_{-n}-z^{-n}| & = 0, \quad |z|>1 .
\label{ssBC:d}
\end{align}
Their linear independence is guaranteed by $ W[\Psi_{-n},\Psi_{+n}]=1 $ for all $ n\geq 0 $.
We have the exact correspondences with Geronimo and Case, $ \phi_{\pm,n}(z)=\phi_{\pm}(Z,n) $ and
$ \hat{\phi}_{\pm,n}(z)=\hat{\phi}_{\pm}(Z,n) $.
\end{definition}

On $ |z|=1 $ we decompose the polynomial solution $ \Psi_{n} $ of (\ref{PsiRecur}) into $ \Psi_{\pm n} $
components
\begin{equation*}
   \Psi_{n} = \begin{pmatrix} \varphi_{n}(z) \cr \bar{\varphi}_{n}(z^{-1}) \end{pmatrix}
  = f_{-}(z)\Psi_{+n}(z)+f_{+}(z)\Psi_{-n}(z) ,
\end{equation*}
introducing the Jost functions $ f_{\pm}(z) $ as the coefficients of this decomposition.
Inverting these relations one finds that
\begin{align*}
   f_{+}(z) & = -\varphi_{n}(z)\hat{\phi}_{+n}(z)+\bar{\varphi}_{n}(z^{-1})\phi_{+n}(z) ,
   \\
   f_{-}(z) & =  \varphi_{n}(z)\phi_{-n}(z)-\bar{\varphi}_{n}(z^{-1})\hat{\phi}_{-n}(z) .
\end{align*}
As a consequence we have the limiting relations, which are a direct analog of Eq. (III.13) in \cite{GC_1979}
\begin{equation*}
   f_{+}(z) = \lim_{n \to \infty}z^n\bar{\varphi}_{n}(z^{-1}) =: \varphi^{*}_{\infty}(z), \quad |z|<1 ,
\end{equation*}
\begin{equation*}
   f_{-}(z) = \lim_{n \to \infty}z^{-n}\varphi_{n}(z), \quad |z|>1 .
\end{equation*}
These limits can be shown to exist because $ f_{+} = W[\Psi_{n},\Psi_{+n}] $ and 
$ f_{-} = -W[\Psi_{n},\Psi_{-n}] $, and the property of the discrete Wronskian given in Proposition \ref{wronskian}.
We note the special values $ f_{+}(0) = f_{-}(\infty) = \kappa_{\infty} $.
In the notations of Wu \cite{Wu_1966} we can identify $ P(z) = f_{+}(z) $ and
$ Q(z^{-1}) = f_{-}(z) $. In the notations of Borodin and Okounkov \cite{BO_2000a} we identify
$ f_{+}(z) = \exp(-[V(-z)]_{>0}) $ and $ f_{-}(z) = \exp(-[V(-z)]_{<0}) $. 

A central orthogonality result of \cite{GC_1979} that generalises without change is the
following.
\begin{corollary}[\cite{GC_1979}]\label{orthog}
Let us assume $ f_{+}(z) \neq 0 $ for $ |z|<1 $ and $ f_{-}(z) \neq 0 $ for $ |z|>1 $.
For $ n \geq m $ the bi-orthogonal polynomials $ \{\varphi_{n},\bar{\varphi}_{n}\}^{\infty}_{n=0} $
satisfy
\begin{equation*}
    \int_{\T}\frac{dz}{2\pi iz} \frac{\varphi_{n}(z)\bar{\varphi}_{m}(z^{-1})}{f_{+}(z)f_{-}(z)} = 0 .
\end{equation*}
Thus we have a factorisation formula, up to a constant, for the weight 
\begin{equation}
    w(z) \propto \frac{1}{f_{+}(z)f_{-}(z)} .
\label{wgt_factor}
\end{equation}
\end{corollary}
\begin{proof}
This follows from the proof given in \cite{GC_1979} with the obvious modifications.
\end{proof}

The scattering solutions defined in Definition \ref{SSdefn} can in fact be related to the associated 
functions previously defined in \textsection \ref{BOPS}.
\begin{proposition}\label{IdentSS}
We can identify the scattering solutions (\ref{scatterS}) together with the boundary
conditions (\ref{ssBC:a}-\ref{ssBC:d}) with the associated functions
\begin{align*}
  \phi_{+n}(z) & = \tfrac{1}{2}f_{+}(z)\epsilon^{<}_{n}(z) ,\quad |z|<1 ,
  \\
  \hat{\phi}_{+n}(z) & = -\tfrac{1}{2}f_{+}(z)z^{-n}\epsilon^{*<}_{n}(z) ,\quad |z|<1 ,
  \\
  \hat{\phi}_{-n}(z) & = -\tfrac{1}{2}f_{-}(z)\epsilon^{>}_{n}(z) ,\quad |z|>1 ,
  \\
  \phi_{-n}(z) & = \tfrac{1}{2}f_{-}(z)z^{-n}\epsilon^{*>}_{n}(z) ,\quad |z|>1 .
\end{align*}
\end{proposition}
\begin{proof}
Given that the combinations on the right-hand sides satisfy the relation (\ref{PsiRecur})
the correspondences follow from a comparison of the boundary conditions (\ref{ops_epsexp:a}-\ref{ops_epsexp:d})
with (\ref{ssBC:a}-\ref{ssBC:d}).
\end{proof}

We will require expansions of the scattering states in their respective domains of
analyticity, with respect to the monomials
\begin{align*}
   \phi_{+n}(z) & = \sum_{n'\geq n}A_{1}(n,n')z^{n'} ,
  \\
   \hat{\phi}_{+n}(z) & = \sum_{n'\geq 1}A_{2}(n,n')z^{n'} ,
  \\
   \phi_{-n}(z) & = \sum_{n'\geq n}\bar{A}_{1}(n,n')z^{-n'} ,
  \\
   \hat{\phi}_{-n}(z) & = \sum_{n'\geq 1}\bar{A}_{2}(n,n')z^{-n'} .
\end{align*}

In conformance with the argument made in \cite{GC_1979} we need to define a scattering function.
\begin{definition}
The scattering function is defined on $ |z|=1 $ as
\begin{equation*}
   S(z) \mathrel{\mathop:}= \frac{f_{-}(z)}{f_{+}(z)}
\end{equation*}
and the sequence of Fourier coefficients $ F_{m}, \bar{F}_{m} $, $ m\in\Z $ by
\begin{equation*}
   F_{m} \mathrel{\mathop:}= \int_{T}\frac{dz}{2\pi iz} z^{m}S(z) ,
   \qquad
   \bar{F}_{m} \mathrel{\mathop:}= \int_{T}\frac{dz}{2\pi iz} z^{-m}\frac{1}{S(z)} .
\end{equation*}
\end{definition}
In comparing with the work of Borodin and Okounkov \cite{BO_2000a} we have $ S(z) = \exp(-V^*(z)) $.
For systems of orthogonal polynomials on the unit circle one has $ |S(z)|=1 $ on $ |z|=1 $ 
and $ \bar{F}_m $ is the complex conjugate of $ F_m $ however this is no longer true for bi-orthogonal polynomial systems.

We can now deduce a system of linear relations for the coefficients defined above which
constitute a discrete analog of the Marchenko integral equations in the inverse scattering theory.
\begin{proposition}\label{MarchenkoEqn}
The coefficients $ A_{1}, \bar{A}_{2}, \bar{A}_{1}, A_{2} $ satisfy the linear relations
\begin{align}
  \frac{\kappa_{n}}{\kappa_{\infty}}\delta_{m,n} & = A_{1}(n,m)+\sum_{n'\geq 1}\bar{A}_{2}(n,n')\bar{F}_{n'+m}, \quad m\geq n\geq 1,
\label{Alinear:a}  \\
  0 & = \bar{A}_{2}(n,m)+\sum_{n'\geq n}A_{1}(n,n')F_{n'+m}, \quad m,n > 0,
\label{Alinear:b}  \\ 
  \frac{\kappa_{n}}{\kappa_{\infty}}\delta_{m,n} & = \bar{A}_{1}(n,m)+\sum_{n'\geq 1}A_{2}(n,n')F_{n'+m}, \quad m\geq n\geq 1,
\label{Alinear:c}  \\
  0 & = A_{2}(n,m)+\sum_{n'\geq n}\bar{A}_{1}(n,n')\bar{F}_{n'+m}, \quad m,n > 0.
\label{Alinear:d}
\end{align}
\end{proposition}
\begin{proof}
These relations follow from a recasting of the jump equations (\ref{jump},\ref{jump*}) using the
factorisation result of Corollary \ref{orthog} and Proposition \ref{IdentSS}. 
Thus to deduce (\ref{Alinear:a}-\ref{Alinear:d}) we employ
\begin{equation}
   \frac{\varphi_{n}(z)}{f_{-}(z)} = \phi_{+n}(z)+\frac{f_{+}(z)}{f_{-}(z)}\hat{\phi}_{-n}(z) ,
\label{factorJ:a}
\end{equation}
\begin{equation}
   \frac{\varphi_{n}(z)}{f_{+}(z)} = \hat{\phi}_{-n}(z)+\frac{f_{-}(z)}{f_{+}(z)}\phi_{+n}(z) ,
\label{factorJ:b}
\end{equation}
\begin{equation}
   \frac{\bar{\varphi}_{n}(z^{-1})}{f_{+}(z)} = \phi_{-n}(z)+\frac{f_{-}(z)}{f_{+}(z)}\hat{\phi}_{+n}(z) ,
\label{factorJ:c}
\end{equation}
\begin{equation}
   \frac{\bar{\varphi}_{n}(z^{-1})}{f_{-}(z)} = \hat{\phi}_{+n}(z)+\frac{f_{+}(z)}{f_{-}(z)}\phi_{-n}(z) ,
\label{factorJ:d}
\end{equation}
respectively.
\end{proof}

We further define kernel matrices $ G_{l,m}, \bar{G}_{l,m} $ (which differ from those given by Eq. (V.11)
of \cite{GC_1979} being the transpose of theirs - this appears to be a typographical error) using the scattering Fourier coefficients
\begin{align*}
   G_{l,m} & \mathrel{\mathop:}= -\sum_{m'\geq 1} \bar{F}_{l+m'}F_{m+m'},
  \\
   \bar{G}_{l,m} & \mathrel{\mathop:}= -\sum_{m'\geq 1} F_{l+m'}\bar{F}_{m+m'} = G_{m,l} .
\end{align*}
We are now in a position to solve the linear system given in the previous proposition.
\begin{proposition}[\cite{GC_1979}]\label{MarchenkoSoln}
The solutions of the discrete Marchenko equations (\ref{Alinear:a}-\ref{Alinear:d}) for the norms and
reflection coefficients are given by
\begin{equation}
   \frac{\kappa_{\infty}^2}{\kappa_{n}^2} = \frac{\det[1+G]^{\infty}_{n+1}}{\det[1+G]^{\infty}_{n}} ,
\label{MEqn:a}
\end{equation}
where
\begin{equation*}
  \det[1+G]^{\infty}_{n} = \det[1+\bar{G}]^{\infty}_{n} 
  = \det
    \begin{pmatrix}
        1+G_{n,n} & G_{n,n+1} & G_{n,n+2} & G_{n,n+3} & \ldots \cr
        G_{n+1,n} & 1+G_{n+1,n+1} & G_{n+1,n+2} & G_{n+1,n+3} & \ldots \cr
        G_{n+2,n} & G_{n+2,n+1} & 1+G_{n+2,n+2} & G_{n+2,n+3} & \ldots \cr
        \vdots & \vdots & \vdots & \vdots & \ddots \cr
    \end{pmatrix} .
\end{equation*}
In addition we have
\begin{equation}
  r_{n+1} = \frac{1}{\det[1+G]^{\infty}_{n+1}} 
  \det
    \begin{pmatrix}
        F_{n+1} & F_{n+2} & F_{n+3} & \ldots \cr
        G_{n+1,n} & 1+G_{n+1,n+1} & G_{n+1,n+2} & \ldots \cr
        G_{n+2,n} & G_{n+2,n+1} & 1+G_{n+2,n+2} & \ldots \cr
        G_{n+3,n} & G_{n+3,n+1} & G_{n+3,n+2} & \ldots \cr
        \vdots & \vdots & \vdots & \ddots \cr
    \end{pmatrix} ,
\label{MEqn:b}
\end{equation}
and
\begin{equation}
  \bar{r}_{n+1} = \frac{1}{\det[1+\bar{G}]^{\infty}_{n+1}} 
  \det
    \begin{pmatrix}
        \bar{F}_{n+1} & \bar{F}_{n+2} & \bar{F}_{n+3} & \ldots \cr
        \bar{G}_{n+1,n} & 1+\bar{G}_{n+1,n+1} & \bar{G}_{n+1,n+2} & \ldots \cr
        \bar{G}_{n+2,n} & \bar{G}_{n+2,n+1} & 1+\bar{G}_{n+2,n+2} & \ldots \cr
        \bar{G}_{n+3,n} & \bar{G}_{n+3,n+1} & \bar{G}_{n+3,n+2} & \ldots \cr
        \vdots & \vdots & \vdots & \ddots \cr
    \end{pmatrix} .
\label{MEqn:c}
\end{equation}
\end{proposition}
\begin{proof}
For convenience we define the coefficient ratios
\begin{equation*}
   a(n,m) = \frac{A_{1}(n,m)}{A_{1}(n,n)}, \qquad \bar{a}(n,m) = \frac{\bar{A}_{1}(n,m)}{\bar{A}_{1}(n,n)} .
\end{equation*}
Summarising (\ref{Alinear:a}-\ref{Alinear:d}) we now have for $ m>n\geq 0 $ the linear relations
\begin{align}
   0 & = a(n,m)+G_{m,n}+\sum^{\infty}_{l=n+1} G_{m,l}a(n,l) ,
\label{newLinear:a}\\
   0 & = \bar{a}(n,m)+\bar{G}_{m,n}+\sum^{\infty}_{l=n+1} \bar{G}_{m,l}\bar{a}(n,l) ,
\label{newLinear:b}
\end{align}
and for $ m=n\geq 0 $
\begin{align}
   \frac{\kappa_{n}}{\kappa_{\infty}A_{1}(n,n)} & = 1+G_{n,n}+\sum^{\infty}_{l=n+1} G_{n,l}a(n,l) ,
\label{newLinear:c}\\
   \frac{\kappa_{n}}{\kappa_{\infty}\bar{A}_{1}(n,n)} & = 1+\bar{G}_{n,n}+\sum^{\infty}_{l=n+1} \bar{G}_{n,l}\bar{a}(n,l) .
\label{newLinear:d}
\end{align}
From the solutions for $ A_{1}, A_{2}, \bar{A}_{1}, \bar{A}_{2} $ one can recover the leading polynomial 
normalisation and reflection coefficients
\begin{gather*}
   \frac{\kappa_{\infty}}{\kappa_{n}} = A_{1}(n,n) = \bar{A}_{1}(n,n) ,
   \\
   r_{n+1} = -\frac{\bar{A}_{2}(n,1)}{\bar{A}_{1}(n,n)} ,
   \qquad
   \bar{r}_{n+1} = -\frac{A_{2}(n,1)}{A_{1}(n,n)} .
\end{gather*}
Using Cramer's rules we can solve the linear relations (\ref{newLinear:a},\ref{newLinear:b}) for 
$ a(n,l), \bar{a}(n,l) $ in terms of determinants and the results are (\ref{MEqn:a},\ref{MEqn:b},\ref{MEqn:c}).
\end{proof}
The Neumann expansion of the basic determinant
\begin{equation*}
    \det[1+G]^{\infty}_{n} \mathop{\sim} 1 + \sum_{n_1\geq n} G_{n_1,n_1}
     + \sum_{n_2>n_1\geq n} \det[G_{n_j,n_k}]_{j,k=1,2} + \dots ,
\end{equation*}
is directly applicable to the formulae given in \textsection \ref{DFD}, (\ref{KV1}) and (\ref{Ddet}),
after noting that $ G_{i,j} = -K(i,j) $. It is this object which specifies the Toeplitz determinant (\ref{Tdet})
\begin{equation*}
   I_{n} = \frac{I_{0}}{\kappa^{2n}_{\infty}}\frac{\det[1+G]^{\infty}_{n}}{\det[1+G]^{\infty}_{0}} .
\end{equation*}

\section{Equality of the Fredholm Determinants}
\label{IM}
\setcounter{equation}{0}

\subsection{$ \lambda$-Extended Riemann-Hilbert problem}
We now pose a new problem which is an extension of the jump conditions given in Proposition \ref{RHjump} by
introducing a parameter $ \lambda $ through the mappings
$ \epsilon^{<}_{n} \mapsto \lambda^2\epsilon^{<}_{n} $ and $ \epsilon^{*<}_{n} \mapsto \lambda^2\epsilon^{*<}_{n} $.
We have two new functional equations where
the associated functions defined in the interior of the unit circle $ \epsilon^{<}_{n}(z), \epsilon^{*<}_{n}(z) $
differ from those in the exterior $ \epsilon^{>}_{n}(z), \epsilon^{*>}_{n}(z) $ on $ |z|=1 $ through
the jump conditions
\begin{gather}
   w(z)\varphi_{n}(z) = -\tfrac{1}{2}\epsilon^{>}_{n}(z)+\lambda^2\tfrac{1}{2}\epsilon^{<}_{n}(z) ,
\label{lambdaJump}\\
   w(z)\varphi^{*}_{n}(z) = \tfrac{1}{2}\epsilon^{*>}_{n}(z)-\lambda^2\tfrac{1}{2}\epsilon^{*<}_{n}(z) .
\label{lambdaJump*}
\end{gather}
Given a weight $ w(z) $
such a system of functional equations can be taken as the defining relations of the bi-orthogonal 
polynomial system and this is essentially the Riemann-Hilbert approach (see Proposition 2.8 of \cite{FW_2006a}). 
The solutions of this extended Riemann-Hilbert problem, $  \varphi_{n}, \varphi^{*}_{n}, \epsilon_{n}, \epsilon^{*}_{n} $,
acquire a $ \lambda $ dependence as do their coefficients and therefore the reflection coefficients.
However we do not address here the question of a $ \lambda $ dependent weight from which one
could construct these coefficients - this will be the subject of our companion work \cite{WF_2011a}.

Having set ourselves the task of solving this new problem we retrace the steps taken in Section \ref{ST} that logically 
follow from using Equations (\ref{lambdaJump}) and (\ref{lambdaJump*}).
We find that Proposition \ref{MarchenkoEqn} now becomes the following.
\begin{proposition}
The $ \lambda $-extended coefficients $ A_{1}, \bar{A}_{2}, \bar{A}_{1}, A_{2} $ satisfy the linear relations
\begin{align*}
  \frac{\kappa_{n}}{\kappa_{\infty}}\delta_{m,n} & = A_{1}(n,m)+\sum_{n'\geq 1}\bar{A}_{2}(n,n')\bar{F}_{n'+m}, \quad m\geq n\geq 1,
\\
  0 & = \bar{A}_{2}(n,m)+\lambda^2\sum_{n'\geq n}A_{1}(n,n')F_{n'+m}, \quad m,n > 0,
\\ 
  \frac{\kappa_{n}}{\kappa_{\infty}}\delta_{m,n} & = \bar{A}_{1}(n,m)+\lambda^2\sum_{n'\geq 1}A_{2}(n,n')F_{n'+m}, \quad m\geq n\geq 1,
\\
  0 & = A_{2}(n,m)+\sum_{n'\geq n}\bar{A}_{1}(n,n')\bar{F}_{n'+m}, \quad m,n > 0.
\end{align*}
\end{proposition}
\begin{proof}
The alteration to the proof of Proposition \ref{MarchenkoEqn} consists of the replacement of 
(\ref{factorJ:b}) and (\ref{factorJ:c}) by the equations
\begin{equation*}
   \frac{\varphi_{n}(z)}{f_{+}(z)} = \hat{\phi}_{-n}(z)+\lambda^2\frac{f_{-}(z)}{f_{+}(z)}\phi_{+n}(z) ,
\end{equation*}
\begin{equation*}
   \frac{\bar{\varphi}_{n}(z^{-1})}{f_{+}(z)} = \phi_{-n}(z)+\lambda^2\frac{f_{-}(z)}{f_{+}(z)}\hat{\phi}_{+n}(z) ,
\end{equation*}
respectively and the same reasoning applied.
\end{proof}

Consequently Proposition \ref{MarchenkoSoln} is extended in the following way.
\begin{proposition}
The solutions of the $ \lambda $-extended discrete Marchenko equations for the norms and
reflection coefficients are given by (\ref{MEqn:a}), (\ref{MEqn:b}), (\ref{MEqn:c}) respectively
with the substitutions $ F_{m}, \bar{F}_{m} \mapsto \lambda F_{m}, \lambda\bar{F}_{m} $ and
$ G_{m,n}, \bar{G}_{m,n} \mapsto \lambda^2G_{m,n}, \lambda^2\bar{G}_{m,n} $.
As a consequence the coefficients have expansions
\begin{gather}
   \frac{\kappa_{\infty}^2}{\kappa_{n}^2} \mathop{\sim}\limits_{\lambda \to 0} 1 - \lambda^2 G_{n,n} + \lambda^4 \sum_{n_1\geq n} G_{n,n_1}G_{n_1,n} + \dots ,
\label{lambdaCft:a}\\
  r_{n+1} \mathop{\sim}\limits_{\lambda \to 0} \lambda F_{n+1} - \lambda^3 \sum_{n_1\geq n+1} F_{n_1+1}G_{n_{1},n} + \ldots ,
\label{lambdaCft:b}\\
  \bar{r}_{n+1} \mathop{\sim}\limits_{\lambda \to 0} \lambda \bar{F}_{n+1} - \lambda^3 \sum_{n_1\geq n+1} \bar{F}_{n_1+1}\bar{G}_{n_{1},n} + \ldots .
\label{lambdaCft:c}
\end{gather}
\end{proposition}
\begin{proof}
For $ m>n\geq 0 $ the linear relations (\ref{newLinear:a}) and (\ref{newLinear:b}) are replaced by
\begin{align*}
   0 & = a(n,m)+\lambda^2 G_{m,n}+\lambda^2\sum^{\infty}_{l=n+1} G_{m,l}a(n,l) ,
\\
   0 & = \bar{a}(n,m)+\lambda^2 \bar{G}_{m,n}+\lambda^2\sum^{\infty}_{l=n+1} \bar{G}_{m,l}\bar{a}(n,l) ,
\end{align*}
respectively and for $ m=n\geq 0 $ Equations (\ref{newLinear:c}) and (\ref{newLinear:d}) are replaced by
\begin{align*}
   \frac{\kappa_{n}}{\kappa_{\infty}A_{1}(n,n)} & = 1+\lambda^2 G_{n,n}+\lambda^2\sum^{\infty}_{l=n+1} G_{n,l}a(n,l) ,
  \\
   \frac{\kappa_{n}}{\kappa_{\infty}\bar{A}_{1}(n,n)} & = 1+\lambda^2 \bar{G}_{n,n}+\lambda^2\sum^{\infty}_{l=n+1} \bar{G}_{n,l}\bar{a}(n,l) ,
\end{align*}
respectively. These are solved in the same way as indicated in the proof of Proposition \ref{MarchenkoSoln}.
\end{proof}
The Neumann expansion of the basic determinant is now
\begin{equation*}
    \det[1+\lambda^2G]^{\infty}_{n} \mathop{\sim}\limits_{\lambda \to 0} 1 + \lambda^2\sum_{n_1\geq n} G_{n_1,n_1}
     + \lambda^4 \sum_{n_2>n_1\geq n} \det[G_{n_j,n_k}]_{j,k=1,2} + \dots ,
\end{equation*}
and this specifies the extended Toeplitz determinant (\ref{Tdet})
\begin{equation*}
   I_{n} = \frac{I_{0}}{\kappa^{2n}_{\infty}}\frac{\det[1+\lambda^2G]^{\infty}_{n}}{\det[1+\lambda^2G]^{\infty}_{0}} .
\end{equation*}

\subsection{Iterative solution of $\lambda$-extended jump conditions}
In the final step of the proof of Proposition \ref{Equality} we return to the expression of the form factors
as multiple integrals given in \textsection \ref{FFE}.
The calculation of the form factors (\ref{ff_Ieven},\ref{ff_Iodd}) by Lyberg and McCoy follows the 
method of Wu \cite{Wu_1966} which is essentially the application of Wiener-Hopf arguments combined 
with an iterative solution to the jump condition (\ref{lambdaJump*}). The last part of our proof formalises 
this assertion.
\begin{corollary}
Consider the low temperature case $ \K>1 $. With
\begin{equation*}
   f_{+}(z) = (1-\K^{-1}z)^{1/2}, \qquad f_{-}(z) = (1-\K^{-1}z^{-1})^{-1/2} ,
\end{equation*}
the solution of the coupled set of equations
\begin{gather}
  \kappa_{n}\varphi^{*}_{n}(z) = f_{+}(z) \left( \left[ f_{-}(z) \right]_{\geq 0}-\lambda^2\tfrac{1}{2}\kappa_{n}\left[ f_{-}(z)\epsilon^{*<}_{n}(z) \right]_{\geq 0} \right) ,
\label{iterate:e}
  \\
  \tfrac{1}{2}\kappa_{n}\epsilon^{*>}_{n}(z)-1 = \frac{1}{f_{-}(z)} 
              \left( -\left[ f_{-}(z) \right]_{< 0}+\lambda^2\tfrac{1}{2}\kappa_{n}\left[ f_{-}(z)\epsilon^{*<}_{n}(z) \right]_{< 0} \right) ,
\label{iterate:f}
  \\
  \kappa_{n}z^n\varphi^{*}_{n}(z^{-1}) 
  = f_{-}(z^{-1}) \left( \left[ z^nf_{+}(z^{-1}) \right]_{\geq 0}+\left[ z^nf_{+}(z^{-1})(\tfrac{1}{2}\kappa_{n}\epsilon^{*>}_{n}(z^{-1})-1) \right]_{\geq 0} \right) ,
  \\
  \tfrac{1}{2}\kappa_{n}z^n\epsilon^{*<}_{n}(z^{-1}) = \frac{1}{f_{+}(z^{-1})} 
              \left( \left[ z^nf_{+}(z^{-1}) \right]_{< 0}+\left[ z^nf_{+}(z^{-1})(\tfrac{1}{2}\kappa_{n}\epsilon^{*>}_{n}(z^{-1})-1) \right]_{< 0} \right) ,
\end{gather}
developed as an expansion in $ \lambda^2 $ gives the form factor expansion formula
(\ref{ffexp2}) with (\ref{ff_Ieven}).
\end{corollary}
\begin{proof}
Employing the usual Wiener-Hopf arguments we can deduce the coupled functional equations
\begin{gather}
  \kappa_{n}\varphi^{*}_{n}(z) = f_{+}(z) \left( \left[ f_{-}(z) \right]_{\geq 0}-\tfrac{1}{2}\kappa_{n}\left[ f_{-}(z)\epsilon^{*<}_{n}(z) \right]_{\geq 0} \right) ,
\label{iterate:a}
  \\
  \tfrac{1}{2}\kappa_{n}\epsilon^{*>}_{n}(z)-1 = \frac{1}{f_{-}(z)} 
              \left( -\left[ f_{-}(z) \right]_{< 0}+\tfrac{1}{2}\kappa_{n}\left[ f_{-}(z)\epsilon^{*<}_{n}(z) \right]_{< 0} \right) ,
\label{iterate:b}
  \\
  \kappa_{n}z^n\varphi^{*}_{n}(z^{-1}) 
  = f_{-}(z^{-1}) \left( \left[ z^nf_{+}(z^{-1}) \right]_{\geq 0}+\left[ z^nf_{+}(z^{-1})(\tfrac{1}{2}\kappa_{n}\epsilon^{*>}_{n}(z^{-1})-1) \right]_{\geq 0} \right) ,
\label{iterate:c}
  \\
  \tfrac{1}{2}\kappa_{n}z^n\epsilon^{*<}_{n}(z^{-1}) = \frac{1}{f_{+}(z^{-1})} 
              \left( \left[ z^nf_{+}(z^{-1}) \right]_{< 0}+\left[ z^nf_{+}(z^{-1})(\tfrac{1}{2}\kappa_{n}\epsilon^{*>}_{n}(z^{-1})-1) \right]_{< 0} \right) ,
\label{iterate:d}
\end{gather}
from (\ref{jump*}) and (\ref{wgt_factor}). These are precisely Equations (2.19a), (2.20a), (2.19b) and 
(2.20b) of \cite{Wu_1966} respectively, when re-written in our notation.
The first, second and fourth equations can be considered as a coupled set of functional
equations for $ \varphi^{*}_{n}(z), \epsilon^{*<}_{n}(z) $ and $ \epsilon^{*>}_{n}(z) $.
We now seek an iterative solution to the second and fourth equations starting with the initial term
involving $ \epsilon^{*<}_{n}(z) $ neglected in (\ref{iterate:b}), i.e. a perturbation expansion with
the parameter $ \lambda^2 $ as a book-keeping parameter recording the order of approximation
in this solution. This means that we replace (\ref{iterate:a}) and (\ref{iterate:b}) with
(\ref{iterate:e}) and (\ref{iterate:f}) respectively,
which is expressed by the mapping $ \epsilon^{*<}_{n}(z) \mapsto \lambda^2\epsilon^{*<}_{n}(z) $.
This is the essential idea behind the asymptotic analysis employed by Wu and also the derivation 
of the form factor expansion by Lyberg and McCoy, who employed Cauchy integral representations for the 
positive and negative parts, and which ultimately led to the formulae (\ref{ff_Ieven}) and (\ref{ff_Iodd}). 
By introducing $ \lambda $ in this way we can trace it through the workings of
\cite{LM_2007} and find that it appears precisely as the correct pre-factor of the terms in the
form factor expansion.
\end{proof}

This completes the proof of Proposition \ref{Equality} and we close this part of our study by making two remarks.
\begin{remark}
The two Fredholm determinant evaluations for the diagonal correlations of the Ising model reported here
don't have direct relevance to the multi-variable integral formulas for the form factors of 
correlations away from the diagonal, such as those in \cite{WMTB_1976} or 
\cite{ONGP_2001}, and more work needs to be done to relate these apparently different representations.
\end{remark}
\begin{remark}
Both our Fredholm determinant formulae generalise beyond the diagonal correlations of the Ising model
and are valid for symbols or weights that are of regular semi-classical type.
Therefore all of our results would apply to the row correlations as well with the obvious modifications.
\end{remark}

\subsection{Explicit examples}
The general theory given above can now be applied to give explicit formulae for the Ising model diagonal 
correlations in the low temperature regime $ \K>1 $, due to the fact that the winding number vanishes
and the conditions on the weight and scattering data ensure the validity of this theory. We will evaluate a number
of initial data explicitly and this allows us to make contact with the results in a companion study
\cite{WF_2011a} where the same quantities have been calculated using the fact that our bi-orthogonal polynomial system
is also an isomonodromic system identified with Picard's solution of the sixth Painlev\'e equation.

According to Corollary \ref{orthog} for $ |\K|>1 $ we factorise the weight (\ref{diagonal:a}) into the interior and 
exterior Jost functions
\begin{equation*}
   f_{+}(z) = (1-\K^{-1}z)^{1/2}, \quad f_{-}(z) = (1-\K^{-1}z^{-1})^{-1/2} .
\end{equation*}
In this example $ \kappa_{\infty}=1 $ as $ c_0=0 $.
A straight-forward calculation reveals that the Fourier coefficients with $ |\K|>1 $ and $ n\in \Z $ are
\begin{align*}
   F_{n} & = \frac{\displaystyle\Gamma(|n|+\tfrac{1}{2})}{\displaystyle\pi^{1/2}|n|!}\K^{-|n|} {}_2F_1(\tfrac{1}{2},|n|+\tfrac{1}{2};|n|+1;\K^{-2}) ,
  \\
 \bar{F}_{n} & = -\frac{\displaystyle\Gamma(|n|-\tfrac{1}{2})}{\displaystyle2\pi^{1/2}|n|!}\K^{-|n|} {}_2F_1(-\tfrac{1}{2},|n|-\tfrac{1}{2};|n|+1;\K^{-2}) .
\end{align*}
All these coefficients have complete elliptic integral evaluations (using the standard definitions given in 
Chapter 19 of \cite{DLMF} or Chapter 13 of \cite{EMOT_II}) as the following examples illustrate
\begin{align}
   F_{0} & = \frac{2}{\pi}\eK(t) ,
\label{F0}
   \\
   F_{1} & = F_{-1} = \frac{2}{\pi}t^{-1/2}\left[ \eK(t)-\eE(t) \right] ,
   \\
   \bar{F}_{0} & = \frac{2}{\pi}\left[ (t-1)\eK(t)+2\eE(t) \right] ,
\label{bF0}
   \\
   \bar{F}_{1} & = \bar{F}_{-1} = -\frac{2}{3\pi}t^{-1/2}\left[ (t-1)\eK(t)+(t+1)\eE(t) \right] ,
\end{align}
The preceding evaluations imply that the kernel matrix $ l,m\in \Z $ is given by
\begin{multline*}
   G_{l,m} = \frac{1}{2\pi}\K^{-l-m} \sum^{\infty}_{n=1}\K^{-2n}\frac{\Gamma(n+l-\tfrac{1}{2})\Gamma(n+m+\tfrac{1}{2})}{(n+l)!(n+m)!}
   \\
   \times {}_2F_1(-\tfrac{1}{2},n+l-\tfrac{1}{2};n+l+1;t){}_2F_1(\tfrac{1}{2},n+m+\tfrac{1}{2};n+m+1;t) ,
\end{multline*}
which is convergent for $ |\K|>1 $ as
\begin{equation*}
   F_{n} \mathop{\sim}\limits_{n \to \pm\infty} \pi^{-1/2}(1-t)^{-1/2}|n|^{-1/2}t^{-|n|/2}, \quad
   \bar{F}_{n} \mathop{\sim}\limits_{n \to \pm\infty} -\tfrac{1}{2}\pi^{-1/2}(1-t)^{1/2}|n|^{-3/2}t^{-|n|/2} .
\end{equation*}
However because of the general summation identity $ l,m\in \Z $, $ l\neq m $, $ |t|<1 $
\begin{multline}
   \sum^{\infty}_{n=1}t^{n}\frac{\Gamma(n+l-\tfrac{1}{2})\Gamma(n+m+\tfrac{1}{2})}{(n+l)!(n+m)!}
   {}_2F_1(-\tfrac{1}{2},n+l-\tfrac{1}{2};n+l+1;t){}_2F_1(\tfrac{1}{2},n+m+\tfrac{1}{2};n+m+1;t)
  \\
  = -2\pi \frac{(-\tfrac{1}{2})_{l+1}(\tfrac{1}{2})_{m+1}}{l!m!}\frac{t}{l-m}
    \left[ \frac{1}{m+1}{}_2F_1(\tfrac{1}{2},l+\tfrac{3}{2};l+1;t){}_2F_1(\tfrac{1}{2},m+\tfrac{1}{2};m+2;t)
    \right.
    \\ \left.
          -\frac{1}{l+1}{}_2F_1(\tfrac{1}{2},m+\tfrac{3}{2};m+1;t){}_2F_1(\tfrac{1}{2},l+\tfrac{1}{2};l+2;t) \right] ,
\label{G_offD}
\end{multline}
we find $ G_{l,m}=-K(l,m) $, as given by (\ref{Dkernel}). This provides an alternative derivation of
the discrete Fredholm determinant, to that given in \cite{BO_2000a}. The diagonal elements can be computed
from the recurrence relation
\begin{equation}
   G_{l,l} = \frac{1}{2\pi}\frac{\Gamma(l+\tfrac{1}{2})\Gamma(l+\tfrac{3}{2})}{\Gamma(l+2)^2}t^{l+1}
             {}_2F_1(-\tfrac{1}{2},l+\tfrac{1}{2};l+2;t){}_2F_1(\tfrac{1}{2},l+\tfrac{3}{2};l+2;t)+G_{l+1,l+1} ,
\label{G_onD}
\end{equation}
and the identity $ G_{0,0}+G_{-1,-1}=-1 $.

Again we find that the kernel 
matrix has complete elliptic integral evaluations as the early examples indicate $ |t|<1 $
\begin{align}
  G_{-1,-1} & = -\frac{1}{2\pi^2}\left[ \pi^2+8\eE(t)\eK(t)+4(t-1)\eK(t)^2 \right] ,
\label{G-1-1} \\
  G_{-1,0} & = -G_{0,-1} = \frac{1}{\pi^2}t^{-1/2}\left[ 2\eE(t)^2-4\eE(t)\eK(t)-2(t-1)\eK(t)^2 \right] ,
  \\
  G_{0,0} & = \frac{1}{2\pi^2}\left[ -\pi^2+4(t-1)\eK(t)^2+8\eK(t)\eE(t) \right] ,
\label{G00} \\
  G_{-1,1} & = -3G_{1,-1} = \frac{2}{\pi^2}t^{-1}\left[ \eE(t)-\eK(t) \right]\left[ (t+1)\eE(t)+(t-1)\eK(t) \right] ,
  \\
  G_{0,1} & = 3G_{1,0} = \frac{1}{\pi^2}t^{-1/2}\left[ -6\eE(t)^2-4(t-2)\eE(t)\eK(t)+2(t-1)\eK(t)^2 \right] ,
  \\
  G_{1,1} & = \frac{1}{6\pi^2}t^{-1}\left[ -3\pi^2 t+4(t-1)(3t-2)\eK(t)^2+8(3t-2)\eK(t)\eE(t)+8(t+1)\eE(t)^2 \right] .
\end{align}

\begin{remark}
An important feature of the kernel matrix, given by (\ref{G_offD}) and (\ref{G_onD}), that emerges is the
appearance of the regularised Gauss hypergeometric function rather than the function itself and the 
consequence is that the matrix exists for negative indices $ l,m $, i.e. for all $ l,m \in \Z $. This 
in turn means that it is possible to give definite meaning to the normalisation coefficients $ \kappa_n $
and the reflection coefficients $ r_n, \bar{r}_n $ for negative indices. This is an important theme in
our companion study \cite{WF_2011a}.   
\end{remark}

These examples furnish a check on our results and a link to exact evaluations of the bi-orthogonal polynomial system
coefficients using the connection with Picard's solution to the sixth Painlev\'e equation, see
\cite{MG_2010}, \cite{WF_2011a}. From these works it is known that 
\begin{equation}
  \frac{1}{\kappa_{0}^2} = \frac{I_{1}}{I_{0}} = \sec(x)\left[ \cn(z,t)\dn(z,t)+\sn(z,t)\eZeta(z,t) \right] ,
\label{I1exact}
\end{equation}
in the low temperature phase where the independent variables are
\begin{equation*}
   z:= \frac{2\eK(t)}{\pi}x ,\quad \lambda = \sin(x) ,
\end{equation*}
and $ \lambda=1 $ corresponds to $ x=\pi/2 $.
The $ \cn(z,t), \dn(z,t), \sn(z,t) $ are the standard Jacobian elliptic functions and Jacobi's elliptic zeta function is
\begin{equation*}
   \eZeta(z,t) := E(z,t)-\frac{2\eE(t)}{\pi}x ,
\end{equation*}
with Jacobi's elliptic epsilon function defined by
\begin{equation*}
   E(z,t) := \eE({\rm am}(z,t),t) ,
\end{equation*}
where $ E(z,t) $ is the incomplete second elliptic integral and the elliptic amplitude function 
$ {\rm am}(z,t) $ is defined through the inversion of the incomplete second elliptic integral
\begin{equation*}
   z := F({\rm am}(z,t),t) .
\end{equation*}
By expanding the exact result (\ref{I1exact}) about $ \lambda=0 $ we find that
\begin{multline*}
   \frac{I_{1}}{I_{0}} = 1+\frac{1}{2 \pi ^2}\left[ \pi^2-4(t-1)\eK(t)^2-8\eE(t)\eK(t) \right]\lambda^2
   \\
   +\frac{1}{24\pi^4}\left[ 9\pi^4-40\pi^2(t-1)\eK(t)^2+16(t^2+2t-3)\eK(t)^4+\eE(t)(-80\pi^2\eK(t)+64(t+1)\eK(t)^3) \right] \lambda^4
   \\
   +{\rm O}(\lambda^6) .
\end{multline*}
Clearly this coincides with (\ref{lambdaCft:a}) for $ n=0 $ and (\ref{G00}) to $ {\rm O}(\lambda^2) $. In
addition it is known from the same works that
\begin{equation}
  \frac{1}{\kappa_{-1}^2} = \frac{I_{0}}{I_{-1}} = \sec(x)\frac{1}{\cn(z,t)\dn(z,t)+\sn(z,t)\eZeta(z,t)} ,
\label{I-1exact}
\end{equation}
and the leading terms in the expansion of this ratio give that
\begin{multline*} 
   \frac{I_{0}}{I_{-1}} = 1+\frac{1}{2\pi^2}\left[ \pi^2+8\eE(t)\eK(t)+4(t-1)\eK(t)^2 \right]\lambda^2
   \\
   +\frac{1}{24\pi^4}\left[ 9\pi^4+40\pi^2(t-1)\eK(t)^2+384\eE(t)^2\eK(t)^2+16(5t^2-14t+9)\eK(t)^4+16\eE(t)(5\pi^2\eK(t)+4(5t-7)\eK(t)^3) \right]\lambda^4
   \\
   +{\rm O}(\lambda^6) ,
\end{multline*}
which coincides with (\ref{lambdaCft:a}) for $ n=-1 $ and (\ref{G-1-1}).

Other important examples are $ r_0, \bar{r}_0 $ in the low temperature phase, which are no longer 
unity for $ \lambda\neq 1 $ but rather \cite{WF_2011a}
\begin{gather}
  r_{0} =  \sn(z,t) ,
\label{r0}\\
  \bar{r}_{0} = \frac{1}{\sn(z,t)}\left[ 1+\cn(z,t)\dn(z,t)+\sn(z,t)\eZeta(z,t) \right]\left[ 1-\cn(z,t)\dn(z,t)-\sn(z,t)\eZeta(z,t) \right] .
\label{Cr0}
\end{gather}
In the limit $ \lambda \to 1^{-}$ we note that 
$ r_{0}, \bar{r}_{0} \to 1 $ using the elliptic functions expansions
\begin{align*}
  \cn(z,t) & = -\frac{2}{\pi}\sqrt{1-t}\;\eK(t)(x-\tfrac{\pi }{2})+{\rm O}(x-\tfrac{\pi}{2})^3 ,
\\
  \dn(z,t) & = \sqrt{1-t}+{\rm O}(x-\tfrac{\pi}{2})^2 ,
\\
  \sn(z,t) & = 1+{\rm O}(x-\tfrac{\pi}{2})^2  ,
\\
  Z(x,t) & = -\frac{2}{\pi}\left[ \eE(t)+(t-1)\;\eK(t) \right](x-\tfrac{\pi}{2})+{\rm O}(x-\tfrac{\pi}{2})^3 .
\end{align*}
By expanding (\ref{r0}) and (\ref{Cr0}) about $ \lambda=0 $ we have
\begin{equation*}
   r_{0} = \frac{2}{\pi}\eK(t)\lambda +\frac{1}{3\pi^3}\eK(t)\left[ \pi^2-4(t+1)\eK(t)^2 \right]\lambda^3+{\rm O}(\lambda^5) ,
\end{equation*}
and
\begin{multline*}
   \bar{r}_{0} = \frac{1}{\pi }\left[ 2(t-1)\eK(t)+4\eE(t) \right]\lambda
   \\
    +\frac{1}{3 \pi ^3}\left[ -24\eE(t)^2\eK(t)+2\eE(t)\left(\pi^2-12(t-1)\eK(t)^2\right)+(t-1)\eK(t)\left(\pi^2-4(t-1)\eK(t)^2\right) \right]\lambda^3
    +{\rm O}(\lambda^5) .
\end{multline*}
These coincide with (\ref{lambdaCft:b}) and (\ref{lambdaCft:c}) respectively to $ {\rm O}(\lambda) $ using 
the evaluations (\ref{F0}) and (\ref{bF0}).

We can also make contact with the form factor expansions with the Appell function kernel given in
\textsection \ref{CFD} and without any loss of generality we examine the case $ n=0 $.
From the works of Orrick et al \cite{ONGP_2001}, Boukraa et al \cite{BHMMOZ_2007} and Mangazeev and Guttmann \cite{MG_2010} 
we have expressions for the non-trivial boundary values $ I_{0} $ for $ \lambda\neq 1 $
\begin{equation}
   I_{0} = (1-t)^{1/4} \left\{
   \begin{matrix}
      \frac{\displaystyle \vartheta_{4}(x|\tau)}{\displaystyle \vartheta_{4}(0|\tau)} & t=\K^{-2}, \quad \K>1
      \medskip\\
      \medskip
      \frac{\displaystyle \vartheta_{3}(0|\tau)\vartheta_{1}(x|\tau)}{\displaystyle \vartheta_{2}(0|\tau)\vartheta_{4}(0|\tau)} & t=\K^{2}, \quad \K<1
   \end{matrix} \right. ,
\label{I0exact}
\end{equation}
where we adopt the standard definitions of the elliptic theta functions (see Chapter 20 of \cite{DLMF} 
and \textsection 13.19 of \cite{EMOT_II}) and the elliptic nome is defined by
\begin{equation*}
  q=e^{i\pi \tau} ,\quad \text{and} \quad \tau = i\frac{\eK(1-t)}{\eK(t)} .
\end{equation*}
Naturally in the case $ \lambda=1 $ we have $ I_{0}=1 $ as one can verify from (\ref{I0exact})
using the identities
\begin{equation*}
   t = \left[\frac{\vartheta_{2}(0|\tau)}{\vartheta_{3}(0|\tau)}\right]^4 , \quad
   1-t = \left[\frac{\vartheta_{4}(0|\tau)}{\vartheta_{3}(0|\tau)}\right]^4 , \quad
   \vartheta_{1}(\tfrac{1}{2}\pi|\tau) = \vartheta_{2}(0|\tau), \quad \vartheta_{4}(\tfrac{1}{2}\pi|\tau) = \vartheta_{3}(0|\tau) .
\end{equation*}
Expanding the expressions (\ref{I0exact}) about $ \lambda=0 $ we have in the low temperature regime
\begin{multline*}
  (1-t)^{-1/4}I_{0} = 1 + \frac{2}{\pi^2}\eK(t)\left[ \eK(t)-\eE(t) \right]\lambda ^2
  \\
    + \frac{16}{\pi^4}\eK(t)\left[ \frac{\pi^2}{24}\eK(t)-\frac{\pi^2}{24}\eE(t)
                                   +\frac{1}{8}\eK(t)^3+\frac{1}{8}\eK(t)\eE(t)^2-\frac{1}{12}t\eK(t)^3-\frac{1}{4}\eK(t)^2\eE(t) \right]\lambda^4
    +{\rm O}(\lambda^6) ,
\end{multline*}
which concurs with the initial terms of (\ref{lowT_CFredholmDet}) with (\ref{NeumannLowT}) and (\ref{lowT_Ckernel}).
In the high temperature regime we have
\begin{multline*}
  (1-t)^{-1/4}I_{0} = \frac{2}{\pi}\eK(t) \Bigg\{ \lambda + \frac{4}{\pi^2}\left[ \frac{\pi^2}{24}-\frac{1}{6}(t-2)\eK(t)^2-\frac{1}{2}\eK(t)\eE(t) \right]\lambda^3
  \\
  + \frac{16}{\pi^4}\left[ \frac{3\pi ^4}{640}-\frac{\pi^2}{16}\eK(t)\eE(t)-\frac{1}{12}(t-2)\eK(t)^2\left(\frac{\pi ^2}{4}-\eK(t)\eE(t)\right)
                                  +\frac{1}{8}\eK(t)^2\eE(t)^2+\frac{1}{120}\left(t^2-6t+6\right)\eK(t)^4 \right]\lambda^5
  \\
         +{\rm O}(\lambda^7) \Bigg\} ,
\end{multline*}
which agrees with the initial terms given by (\ref{highT_CFredholm}) with (\ref{highT_Ckernel:a},
\ref{highT_Ckernel:b}, \ref{highT_Ckernel:c}).
The form factor expansion coefficients of this correlation function have been given in
\textsection 4 and \textsection 5.2 of \cite{ONGP_2001}, and we observe agreement with the leading
terms given above. In addition these coefficients have been given in Appendix A of \cite{MG_2010}.

\section{Acknowledgments}
This research was supported by the Australian Research Council under grant ID DP0988652. The authors
would like to express their appreciation for the many questions posed and insights shared by Barry
McCoy, Jacques Perk and Jean-Marie Maillard.

\bibliographystyle{plain}
\bibliography{moment,random_matrices,nonlinear}

\end{document}